\documentstyle[12pt]{article}

\let\Bbb\relax
\newfont{\Bb }{msbm10 scaled 1000}
\newfont{\Bbb}{msbm10 scaled 1200}
\newfam\eufam%
\font\euzw=eufm10 scaled 1200%
\font\euac=eufm9%
\textfont\eufam=\euzw\scriptfont\eufam=\euac%
\def\fr{\fam\eufam\euzw}%
\newcommand{\Z}{{\Bbb Z}}
\newcommand{\R}{{\Bbb R}}
\newcommand{\C}{{\Bbb C}}
\newcommand{\Ha}{{\Bbb H}}
\newtheorem{thm}{Theorem}
\newtheorem{prop}{Proposition}
\newtheorem{cor}{Corollary}
\newtheorem{lemma}{Lemma}
\newtheorem{dof}{Definition}
\begin{document}
\begin{center}
{\LARGE Classification of stationary compact\\[0.3cm]
homogeneous special pseudo-K\"ahler manifolds\\[0.3cm]
of semisimple groups}
\end{center}

\vskip 1.0 true cm 
\begin{center}
{\large D.V.\ Alekseevsky$^{\star}$\footnote{Supported by 
Erwin Schr\"odinger Institut (Vienna). e-mail: daleksee@esi.ac.at},
V.\ Cort\'es$^{\star \star}$\footnote{Supported by SFB 256 (Bonn University). e-mail: 
vicente@math.uni-bonn.de}}\end{center}
\noindent
{\small
$^{\star}$ Sophus Lie Center,
Gen.\ Antonova 2 - 99,
117279 Moscow,
Russia\\
$^{\star \star}$ Mathematisches Institut der Universit\"at Bonn,
Beringstr. 1,
53115 Bonn, Germany\\
}

\begin{abstract} 
The variation of the Hodge structure of a Calabi-Yau 3-fold induces
a ca\-no\-ni\-cal K\"ahler metric on its Kuranishi moduli space, 
known as  the Weil-Petersson metric.  Similarly, {\it special} 
pseudo-K\"ahler manifolds correspond to certain (abstract) variations 
of Hodge structure which generalize the above example.
We give the classification of homogeneous special pseudo-K\"ahler 
manifolds of semisimple group with  
compact stabilizer. 

\medskip\noindent
{\it Key words}: Homogeneous K\"ahler manifolds, Special Geometry,
Calabi-Yau 3-folds/Mirror Symmetry

\medskip\noindent   
{\it MSC}: 53C30 
   \end{abstract} 

\section*{Introduction}
The Kuranishi moduli space of a Calabi-Yau threefold carries a canonical
negative definite K\"ahler metric, known as the Weil-Petersson metric.
More generally, any variation of principally polarized Hodge structure of
weight 3 with $h^{3,0} = 1$ gives rise to a canonical negative definite
K\"ahler metric $g$ on the complex manifold $M$ pa\-ra\-me\-tri\-zing
the variation.
The pair $(M,g)$ is considered as {\em formal moduli space}. Generalizing
the construction to the case of arbitrary polarizations one obtains
pseudo-K\"ahler manifolds $(M,g)$ of arbitrary signature. Such manifolds
are called {\em special} pseudo-K\"ahler manifolds. The notion of
special (pseudo-) K\"ahler geometry was introduced by physicists, see
e.g.\ \cite{10}. In physical terminology, special K\"ahler geometry
is the geometry allowed for the coupling of $N = 2$ supergravity to
vector multiplets in $D = 4$ spacetime dimensions. {\em Homogeneous}
special K\"ahler manifolds were first considered by Cecotti \cite{6},
who related them to homogeneous quaternionic K\"ahler manifolds, cf.\
\cite{1}, \cite{11}, \cite{12}, \cite{7}, \cite{3},
\cite{8} and \cite{9}.  The par\-ti\-cu\-lar\-ly interesting 
homogeneous special
K\"ahler manifolds associated to homogeneous real affine cubic hypersurfaces
were completely classified in \cite{11} and \cite{8}. In general,
none of the two factors in the Levi-Malcev decomposition of the full
isometry group of these homogeneous K\"ahler manifolds is trivial,
see \cite{12} and \cite{3}. The reader who has learned about 
special K\"ahler manifolds from \cite{13} may be surprised
to read there exist homogeneous special K\"ahler manifolds. In fact
it follows from \cite{18} that any special K\"ahler manifold {\it in the
sense of} \cite{13} which has a transitive isometry group is necessarily
flat. However this result does not apply to special K\"ahler manifolds
in our sense, which in the terminology of \cite{13} are not   
special K\"ahler manifolds but {\it projective} special K\"ahler manifolds. 

In the following text, we give the classification of homogeneous 
stationary compact special pseudo-K\"ahler manifolds $M$ of a real semisimple
group $G^0$. (A homogeneous manifold $M \cong G^0/H^0$ is called 
stationary compact if its stationary {\it or isotropy} group $H^0$ is compact.)  
This includes the classification of all homogeneous special
K\"ahler manifolds of a real semisimple group. 
The precise definitions of special pseudo-K\"ahler manifolds, formal
moduli spaces etc.\ are given in the first section. In the second 
section, we prove that if $M$ is a stationary compact homogeneous special 
K\"ahler manifold of a real semisimple group $G^0$ then there exists 
a symplectic module $V$ of the complexified group
$G = (G^0)^{\mbox{\C}} \subset GL(V)$ such that the highest weight
vector orbit ${\cal C} \subset V$ is a Lagrangian cone and $M$ is embedded
into its projectivization $P({\cal C}) \subset P(V)$ as an open
$G^0$-orbit. Then, in section 3, we classify all Lagrangian highest
weight vector orbits of complex semisimple groups $G$.
In the last section, we study the real forms of the symplectic $G$-modules
$V$ occurring in section 3, obtaining the solution to our classification
problem.

 In particular, we establish a natural 1-1 correspondence
between  simple complex Lie algebras $\fr l$ different from
${\fr c}_n = {\fr sp}_n(\mbox{\C})$ and homogeneous special K\"ahler
manifolds $G^0/H^0$ of a semisimple Lie group $G^0$ which are
formal moduli spaces. All such formal moduli spaces are
Hermitian symmetric manifolds of nonpositive curvature and 
belong to the following list :
\begin{enumerate}
\item[A)] ${\fr l} ={\fr {sl}}_{n+2}(\mbox{\C}), \qquad
M = \mbox{\C}H^{n-1} = SU_{1,n-1}/S(U_1 \cdot U_{n-1}) $
\item[BD)] ${\fr l} = {\fr {so}}_{n+4}(\mbox{\C}), \qquad
M = (SL_2(\mbox{\R})/SO_2) \times (SO_{2,n-2}/SO_2\cdot SO_{n-2}) $
\item[G)] ${\fr l} ={\fr g}_2(\mbox{\C}), \qquad
M =\mbox{\C}H^1 = SL_2(\mbox{\R})/SO_2$
\item[F)] ${\fr l} ={\fr f}_4(\mbox{\C}), \qquad
M = Sp_3(\mbox{\R})/U_3$
\item[E6)] ${\fr l} ={\fr e}_6(\mbox{\C}), \qquad
M = SU_{3,3}/S(U_3 \cdot U_3)$
\item[E7)] ${\fr l} ={\fr e}_7(\mbox{\C}), \qquad
M = SO^*_{12}/U_6$
\item[E8)] ${\fr l} ={\fr e}_8(\mbox{\C}), \qquad
M = E_7^{(-25)}/E_6 \cdot SO_2.$  
\end{enumerate}

In order to explain the correspondence we recall Wolf's construction
of the symmetric quaternionic K\"ahler manifolds of compact type \cite{25}. 
Let $\fr l$ be a simple complex Lie algebra, $L^0$ the
corresponding compact Lie group with trivial centre, $Sp(1)
\subset L^0$ the 3-dimensional regular subgroup associated
to the highest root and $K = Z_{L^0}(Sp(1))$ its centralizer. 
The homogeneous manifold $W({\fr l}) := L^0/Sp(1) K$ is called
the {\it Wolf space} associated to $\fr l$. The Wolf spaces
are precisely the symmetric quaternionic K\"ahler manifolds
of positive scalar curvature. Next we consider the isotropy
representation $Sp(1)K \rightarrow GL(V)$, $V = T_{[e]}W$, of the Wolf space 
$W = W({\fr l})$ and denote by $\rho = \rho_{\fr l}: K \rightarrow
GL(V)$ its restriction to the subgroup $K$. The $K$-module $V$
carries a $K$-invariant structure of a complex symplectic vector space.
In other words, to any complex simple Lie algebra $\fr l$ 
we have associated a compact Lie group $K$ and a complex
representation $\rho_{\fr l}: K \rightarrow Sp_{\mbox{\C}}(V)$ of symplectic
type. If ${\fr l} \neq {\fr sp}_{n+1}({\mbox{\C}})$ then $V$ admits a 
Lagrangian highest weight vector orbit ${\cal C} = 
\rho_{\fr l}(K)^{\mbox{\C}}v = \mbox{\C}^{\ast}\rho_{\fr l}(K)v = 
\mbox{\R}^{\ast}\rho_{\fr l}(K)v \subset V$. 
(In the case ${\fr l} = {\fr sp}_{n+1}({\mbox{\C}})$ we have that
$V = \mbox{\Ha}^n = \mbox{\C}^{2n}$ is the standard module of 
$K = Sp(n)$, which has an open highest weight vector orbit ${\cal C} = 
\rho_{\fr l}(K)^{\mbox{\C}}v = V -\{ 0\}$.)
The projectivization $P({\cal C})$ of the complex cone $\cal C$ 
is a compact Hermitian symmetric manifold and the dual noncompact
Hermitian symmetric manifold is a formal moduli space. 
Conversely, any homogeneous formal moduli space can be obtained
by this construction. 
 
  It is natural to conjecture that there exist  Calabi-Yau 3-folds 
whose Kuranishi moduli spaces are  the   formal moduli
spaces from the above list. The recent re\-mar\-ka\-ble results of
Looijenga and Lunts \cite{17} concerning the action of a semisimple
Lie group on the cohomology of a K\"ahler manifold  support
this conjecture (due to Mirror Symmetry).  

Finally, we point out a remarkable coincidence with the irreducible 
holonomy groups of torsion free connections of symplectic type, see
\cite{19} and \cite{23}. Let ${\fr l} \neq {\fr sl}_n(\mbox{\C})$
be a complex simple Lie algebra. Then  the complex symplectic
representation $\rho : K \rightarrow Sp_{\mbox{\C}}(V)$ constructed
above is {\it irreducible} and there exists a complex torsionfree
connection with holonomy group $\rho_{\fr l}(K)^{\mbox{\C}} \subset
Sp_{\mbox{\C}}(V) \cong Sp_m(\mbox{\C}) = Sp(\mbox{\C}^{2m})$. Conversely, if 
$Hol \subset Sp_m(\mbox{\C})$, $m>2$, is the holonomy group
of a complex torsionfree connection and if $Hol$ is irreducible
then $Hol = \rho_{\fr l}(K)^{\mbox{\C}}$ for some complex simple
Lie algebra ${\fr l} \neq {\fr sl}_n(\mbox{\C})$.  Notice that
the generic holonomy group $Hol = Sp_m(\mbox{\C})$ corresponds
to the simple Lie algebra $l = {\fr sp}_{m+1}(\mbox{\C})$. 

We thank Simon Salamon for his useful comments on our paper. 
\section{Definitions} Let $(V,\omega )$ be a complex
symplectic  vector space and $\tau :  V \rightarrow V$ a
compatible real structure on $V$, i.e.\ a {\C}-antilinear
involution such that $\omega (\tau x, \tau y) = \overline {\omega (x,y)}$
for all $x,y \in V$. In other words, $(V,\omega )$ is the complexification of
the real symplectic vector space $(V^{\tau}, \omega | V^{\tau} )$;
where $V^{\tau} = \{ x\in V| \tau x = x\}$ is the fixed point set of $\tau$.
Choosing a real symplectic basis for $V$ we can identify $V =
T^{\ast}\mbox{\C}^n$, $\omega = \sum_{i=1}^n dq^i\wedge dp_i$ in
canonical coordinates $(q^i,p_i)$ and $V^{\tau} = T^{\ast}\mbox{\R}^n$.
Note that  the automorphism group of the triple $(V,\omega , \tau )$ is
$$ Aut (V,\omega , \tau ) = \{ \varphi \in GL(V) | \varphi^{\ast} \omega =
\omega \, ,\quad \varphi^{\ast} \tau =
\tau \} \cong Sp_n (\mbox{\R}) \subset GL_{2n} (\mbox{\R})\, .$$
Combining $\omega$ and $\tau$ we define a Hermitian form $\gamma$ on $V$
by:
$$ \gamma (x,y) := \sqrt{-1} \omega (x,\tau y)\, ,\quad x,y \in V\, .$$
In canonical coordinates it has the following expression:
\[ \gamma = \sqrt{-1} \sum_{i=1}^n (dq^i\otimes d\bar{p}_i -
dp_i \otimes d\bar{q}^i) = \sum_{i=1}^n (e^i_- \otimes \bar{e}^i_- -
e^i_+\otimes \bar{e}^i_+)\, ,\]
where $e^i_{\pm} := \frac{1}{\sqrt{2}}(dq^i \pm \sqrt{-1}dp_i)$. In particular,
$\gamma$ has signature $(n,n)$.

\begin{dof}\label{LagrangianDef} Let $(V,\omega , \tau )$
be as above and $\iota : {\cal C} \hookrightarrow V$ a holomorphic immersion
of a connected complex manifold $\cal C$. The pair $({\cal C}, \iota )$
is called
\begin{enumerate}
\item[i)] {\bf conic} if $\lambda \cdot \iota ({\cal C}) \subset \iota ({\cal C})$
for all $\lambda \in \mbox{\C}^{\ast}$,
\item[ii)] {\bf totally isotropic} if $\iota^{\ast} \omega = 0$,
\item[iii)] {\bf Lagrangian} if $({\cal C},\iota )$ is totally isotropic and
$2\dim {\cal C} = \dim V$,
\item[iv)] {\bf pseudo-K\"ahler} if $\iota^{\ast} \gamma$ is nondegenerate,
\item[v)] a {\bf special cone} if conditions i)-iv) hold and
$\gamma (u,u) \neq 0$ for all $u \in \iota ({\cal C})$.
\end{enumerate}
A connected complex submanifold ${\cal C} \subset V$ is called
conic (resp.\ totally isotropic, Lagrangian etc.) if $({\cal C}, \iota )$
is conic (resp.\ totally isotropic, Lagrangian etc.), where $\iota$ is the
canonical inclusion ${\cal C} \hookrightarrow V$.
\end{dof}

Note that $n = \dim V/ 2$ is the maximal possible dimension of a totally
isotropic complex submanifold ${\cal C} \subset V$. To any
conic complex submanifold ${\cal C} \subset V$ we associate its
projectivization $P({\cal C}) \subset P(V)$, which is a complex
submanifold of the projective space $P(V)$. If $\cal C$ is a special
cone, a canonical pseudo-K\"ahler metric $g$ on $M =P({\cal C})$ is
defined as follows:
\[ g_{\pi u} (d\pi \, v, d\pi \, v) := \frac{\gamma (v,v)}{\gamma (u,u)} -
{\left| \frac{\gamma (v,u)}{\gamma (u,u)}\right|}^2\, ,\quad
u\in {\cal C}\, ,v\in T_u{\cal C}\, ,\]
where $\pi : {\cal C} \rightarrow P({\cal C})$ is the natural projection.

\begin{dof} \label{specialDef} The pseudo-K\"ahler manifold
$(M = P({\cal C}),g)$ canonically
associated to a special cone ${\cal C} \subset V$ is called a {\bf special
pseudo-K\"ahler manifold}. $(M,g)$ is called {\bf special  K\"ahler manifold}
if $g$ or $-g$ is positive definite. If $g$ is negative definite, we call
$(M,g)$ a {\bf formal moduli space} (of Calabi-Yau 3-folds).
\end{dof}

\noindent
{\bf Remark 1:} Let $X$ be a (general) {\bf Calabi-Yau 3-fold}, i.e.\
a compact K\"ahler ma\-ni\-fold with holonomy group $SU_3$, $M_X$ its
Kuranishi moduli space and $g_{WP}$ the Weil-Petersson metric on $M_X$.
Then $(M_X,g_{WP})$ is a formal moduli space in the sense of
Def.\ \ref{specialDef}, cf.\ \cite{9}.

\begin{dof}  \label{homogDef} A special pseudo-K\"ahler manifold $(M,g)
\subset P(V)$
is called {\bf homogeneous} if there exists a closed subgroup
$G \subset Aut (V,\omega , \tau ) \cong Sp_n (\mbox{\R})
\subset GL_{2n} (\mbox{\R})$ such that $M = Gp$ is the $G$-orbit of
some point $p \in P(V)$.
\end{dof}

\noindent
Note that, under the assumptions of Def.\ \ref{homogDef}, $G$ acts
transitively on the pseudo-K\"ahler manifold $(M,g)$ by holomorphic
isometries and we can identify $M$ with the coset space $G/G_p$,
$G_p = \{ \varphi \in G|\varphi p = p\}$.

\begin{prop} Let ${\cal C} \subset V$ be a special cone which is not
contained in a linear Lagrangian subspace of $V$ and 
$G \subset Aut (V,\omega , \tau )$ a subgroup preserving $\cal C$.
Then $PG = G/G\cap \{ \pm {id}_V\}$  acts
effectively on $M = P({\cal C})$ and hence we have
an inclusion
 $$PG \subset Aut (M) = \{ \varphi : M\rightarrow M|
\varphi  \: \mbox{is a  holomorphic isometry} \} .$$
\end{prop}

\noindent
{\bf Proof:} We prove that if $\varphi \in Aut(V,\omega , \tau )$ 
preserves $\cal C$ and acts as identity on $M = P({\cal C})$
then $\varphi = \pm {id}_V$. First we show that $\varphi$ is completely
determined by its value at a single point $u \in {\cal C}$.  
Since $\varphi$ acts as identity on $M$ we know that 
$\varphi (u) = \lambda u$ for some $\lambda \in \mbox{\C}^{\ast}$. 
For the same reason the differential $d\varphi_u : T_u{\cal C} \rightarrow 
T_{\lambda u}{\cal C} = T_u{\cal C}$ induces the identity map on
$T_pM$, $p = {\mbox{\C}}^{\ast}u$, and hence $d\varphi_u (v) = \lambda v$ 
for all $v\in T_u{\cal C} \subset V$. This shows that $\varphi (v) = 
d\varphi_u(v) = \lambda v$ for all $v\in T_u{\cal C}$. The nondegeneracy
of $\gamma|T_u{\cal C}$ is equivalent to
$V = T_u {\cal C} \oplus \tau T_u{\cal C}$. Now using $\varphi^{\ast} \tau =
\tau$ we can conclude that $\varphi \in Aut (V,\omega , \tau )$ is
completely determined by its value $\varphi (u) = \lambda u$
at the point $u\in {\cal C}$. In fact, $\varphi$ preserves
the decomposition  $V = T_u {\cal C} \oplus \tau T_u{\cal C}$ and 
acts as scalar multiplication by $\lambda$ and $\bar{\lambda}$ on the 
Lagrangian subspaces $T_u {\cal C}$ and $\tau T_u{\cal C}$ respectively. 
From the fact that $\cal C$ consists of eigenvectors of $\varphi$ it follows 
that $\lambda = \bar{\lambda}$ or ${\cal C} \subset T_u {\cal C}$ or
${\cal C} \subset \tau T_u {\cal C}$. By assumption $\cal C$ is not
contained in a Lagragian subspace. So $\lambda = \bar{\lambda}$ and 
from $\varphi^{\ast}\gamma = \gamma$ it follows that 
$\lambda = \pm 1$. This shows that $\varphi = \pm {id}_V$ 
proving the proposition. 

\begin{cor} \label{closedCor} Let $M = P({\cal C})$ be a homogeneous 
special K\"ahler manifold of the group
$G$ such that $\cal C$ is not contained in a linear Lagrangian subspace. 
Then $PG \subset Aut (M)$ is a closed subgroup with compact 
stabilizer $PG_p$, $p\in M$.
\end{cor}

\noindent
{\bf Proof:} Recall that for any Hermitian manifold $M$ the subgroup
$Aut (M) \subset Isom (M)$ is closed with respect to the compact open topology.
Since $G$ acts transitively on $M$, we can express the closure of 
$PG$ in $Aut (M)$ as $\overline{PG} = PG\overline{PG}_p$. This shows that it is
sufficient to prove the compactness of $G_p$. Consider the compact
group
\[ K_p = \{ \varphi \in Aut (V,\omega , \tau )| \varphi T_u{\cal C}
\subset   T_u{\cal C}, \quad \varphi (u) \in \mbox{\C} u\}
\cong U_{n-1}\times U_1\, ,\]
where $p = \pi u$, $\pi : {\cal C} \rightarrow M = P({\cal C})$ the natural
projection. The claim follows from the fact that $G_p$ is a closed
subgroup of $K_p \subset Aut (V,\omega , \tau )$. This proves the 
corollary. 

Our aim is to classify all homogeneous special pseudo-K\"ahler manifolds
$(M,g)$ of (real) semisimple group $G$ with compact stabilizer $G_p$.
By Corollary \ref{closedCor} this includes the classification of homogeneous
special K\"ahler manifolds of semisimple group.

\section{Compactification of homogeneous
pseu\-do-\-K\"ah\-ler manifolds}
Let $(M,g)$ be a homogeneous special pseudo-K\"ahler manifold
of the real semisimple group $G$ with compact stabilizer $G_p$, $p\in M$.
We consider the complexified group $G^{\mbox{\C}} \subset GL(V)$, i.e.\
the (algebraic) linear group with Lie algebra ${\fr g}^{\mbox{\C}} =
{\fr g} \otimes \mbox{\C} \subset {\fr gl} (V)$, where ${\fr g} =
Lie\, G \subset {\fr gl} (V^{\tau}) \cong {\fr gl}_{2n}(\mbox{\R})$.
The orbit $\overline{M} = G^{\mbox{\C}}p$ of any point $p
\in M$ is a complex submanifold $\overline{M} \subset P(V)$
containing $M$ as an open subset. We will prove that  $\overline{M}$
is precisely the Zariski closure of $M$ and hence a projective homogeneous
algebraic variety.

More generally, let $G^0$ be a connected real semisimple group
acting effectively and transitively on a pseudo-K\"ahler manifold
$(M,g,J)$, $J$ the complex structure, by holomorphic isometries.
We assume that the stabilizer $H^0$ of a point $p\in M$ is compact.

\begin{thm}  \label{redThmI} Under the conditions above any choice of
base point $p\in M$
defines a canonical open holomorphic embedding $M\hookrightarrow G/P$
into a flag manifold of the connected complex semisimple group $G 
= Aut_e  ({\fr g})$, ${\fr g} = ({\fr g}^0) \otimes \mbox{\C}$,
${\fr g}^0 = Lie \, G^0$. Moreover, the stabilizer $H^0 = Z_{G^0}(h)$
is the centralizer of some element $h\in {\fr h}^0 = Lie \, H^0$
and is connected; $M \cong G^0/H^0$ is simply connected.
\end{thm}

\noindent
{\bf Proof:} The (pseudo-) K\"ahler form on $M$ is a $G^0$-invariant
symplectic structure. It follows from Kirillov,
Kostant and Souriau's description of homogeneous symplectic manifolds
that the stabilizer $H^0 = Z_{G^0}(h)$ for some
$h \in {\fr h}^0 = Lie \, H^0$. The maximal compact subgroup $K^0$ of
$G^0$ is connected (see \cite{14} Ch.\ VI, Theorem 2.1) and hence
$H^0 = Z_{K^0}(h)$ is the centralizer of a torus in a compact connected
Lie group. This shows that $H^0$ is connected and $M$ simply connected
(see \cite{14} Ch.\ VII Corollary 2.7).  Moreover, since $G^0$ acts effectively,
we can identify $G^0$ with its adjoint form $Ad \, G^0$.

Any parabolic
subalgebra ${\fr p} \subset
{\fr g} = {\fr g}^0 \otimes \mbox{\C}$ defines a flag manifold $G/P$, where
$G =  Aut_e  ({\fr g})$ is the adjoint form of $\fr g$ and $P$ is its parabolic subgroup
with $Lie \, P = {\fr p}$. If  $P\cap G^0 = H^0$ then we have an open
embedding $M = G^0/H^0 \hookrightarrow G/P$.

The effective action of $G^0$ on $M$ gives rise to a monomorphism
$\rho^0$ of ${\fr g}^0$ into the Lie algebra ${\fr X}_p(M)$ of germs of vector
fields at $p\in M$. Note that ${\fr h}^0 = \{ x\in {\fr g}^0|
\rho^0 (x) (p) = 0\}$. The infinitesimal action $\rho^0:
{\fr g}^0\rightarrow  {\fr X}_p(M)$ by (germs of) holomorphic Killing vector
fields extends naturally to an infinitesimal action $\rho : {\fr g}
\rightarrow  {\fr X}_p(M)$ by holomorphic vector fields:
$\rho (x + iy ) = \rho^0 (x) + J \rho^0 (y)$, $x,y \in {\fr g}^0$.

\begin{lemma} $\rho : {\fr g} \rightarrow {\fr X}_p(M)$ is injective.
\end{lemma}

\noindent
{\bf Proof:} $\ker \rho$ is an ideal and hence a semisimple direct summand
of $\fr g$. It follows that $\ker \rho$ is the complexification of the ideal
${\fr g}^0\cap \ker \rho \subset {\fr g}^0$. Now the lemma follows from
${\fr g}^0\cap \ker \rho = \ker \rho^0 = 0.$  

The next proposition finishes the proof of Theorem \ref{redThmI}.

\begin{prop} \label{basicProp} The infinitesimal stabilizer ${\fr g}_p = \{
x\in {\fr g}| \rho (x) (p) = 0\}$ is a parabolic subalgebra
${\fr g}_p = {\fr p} \subset {\fr g}$ with maximally reductive
subalgebra ${\fr h} = {\fr z}_{{\fr g}}(h)$. The cor\-res\-pon\-ding
parabolic subgroup $P \subset G$ satisfies $P\cap G^0 = H^0$.
\end{prop}

\noindent
{\bf Proof:} The Lie group $H^0 = Z_{G^0}(h)$ contains a compact Cartan
subgroup $A^0$ of $G^0$ such that $h \in {\fr a}^0 = Lie \, A^0 
\subset {\fr h}^0$. Let $R$ denote the root system of $\fr g$ 
with respect to the
Cartan subalgebra ${\fr a} = {\fr a}^0 + i {\fr a}^0 \subset {\fr g}$.
Note that ${\fr a}^0 = \{ x\in {\fr a}| \alpha (x) \in i\mbox{\R}$ for all
$\alpha \in R\}$. Since ${\fr a} \subset {\fr g}_p$, we have
\[ {\fr g}_p = {\fr a} + \sum_{\alpha \in R_p} {\fr g}_{\alpha}\]
for some root subsystem $R_p \subset R$. More precisely,
since ${\fr g}_p \cap {\fr g}^0 = {\fr h}_0$, we have
\[ {\fr g}_p = {\fr h} + {\fr m}\, ,\]
\[{\fr h} = {\fr h}^0 \otimes \mbox{\C} = {\fr z}_{\fr g} (h) =
{\fr a} + \sum_{\alpha \in R_0} {\fr g}_{\alpha}\, , \quad {\fr m} =
\sum_{\alpha \in R_m} {\fr g}_{\alpha}\, ,\]
where $R_0 = \{ \alpha \in R| \alpha (h) = 0\}\subset R_p$ and
$R_m = R_p - R_0$.

\begin{lemma} $R_p = R_0 \stackrel{\cdot}{\cup} R_m$ is the
decomposition of the root system $R_p$ into its maximally symmetric root
subsystem $R_0$ and the asymmetric root subsystem $R_m$;
${\fr g}_p = {\fr h} + {\fr m}$ is the decomposition of the Lie algebra
${\fr g}_p$ into the maximally reductive subalgebra ${\fr h}$ and the nilpotent
ideal $\fr m$.
\end{lemma}

\noindent
{\bf Proof:} Let $\tau$ denote the {\C}-antilinear involutive automorphism
of $\fr g$ with fixed point set ${\fr g}^{\tau} = {\fr g}^0$. The Cartan
subalgebra ${\fr a} = {\fr a}^0 + i{\fr a}^0$, ${\fr a}^0 ={\fr a}^{\tau}$,
is $\tau$-invariant and a straightforward computation shows that
$\tau {\fr g}_{\alpha} = {\fr g}_{{\alpha}^{\tau}}$, $\alpha \in R$,
where ${\alpha}^{\tau}(a) = \overline{\alpha (\tau (a))}$, $a\in {\fr a}$.
Using the fact that $\alpha$ is purely imaginary on ${\fr a}^0$ we
conclude that ${\alpha}^{\tau} = - \alpha$ and hence
$\tau {\fr g}_{\alpha} = {\fr g}_{-\alpha}$. Let $\alpha \in R_p$ and
assume that also $-\alpha \in R_p$, i.e.\
${\fr g}_{\alpha} + {\fr g}_{-\alpha} \subset {\fr p}$. If ${\fr g}_{\alpha} =
\mbox{\C} e_{\alpha}$, then it follows that
\[ 0 \neq \mbox{\R} (e_{\alpha} + \tau e_{\alpha})
= ({\fr g}_{\alpha} + {\fr g}_{-\alpha})^{\tau} \subset
{\fr g}_p\cap {\fr g}^0 = {\fr h}^0 \]
and hence ${\fr g}_{\alpha} + {\fr g}_{-\alpha} \subset {\fr h}$.
This implies that $\alpha \in R_0$. We have proved that $R_0 \subset
R_p$ is the (unique) maximally symmetric root subsystem and that the
complement  $R_m = R_p - R_0$ is asymmetric. {}From the Levi-Malcev
decomposition it follows that there exists a nilpotent ideal ${\fr n}
\subset {\fr g}_p$ complementary to $\fr h$. Since $\fr h$ contains the Cartan
subalgebra $\fr a$, the complement ${\fr m} = \sum_{\alpha \in R_m}
{\fr g}_{\alpha}$ of $\fr h$ in ${\fr g}_p$ is the only $\fr h$-invariant
complement and hence ${\fr m} = {\fr n}$. This shows, in particular,
that $\fr m$ is a subalgebra and hence $R_m$ is a root system
finishing the proof of the lemma. 

We continue the proof of Proposition \ref{basicProp}. The solvable Lie algebra
${\fr a} + {\fr m} \subset {\fr g}$ is contained in some Borel subalgebra
and so $R_m\subset R^+$ for some system of positive roots $R^+ \subset R$.
Put $R_0^+ := R_0 \cap R^+$. Then $R_0^+$ is a system of positive roots
for the root system $R_0$.

\begin{lemma} \label{parabLemma} $R^+ = R_0^+ \stackrel{\cdot}{\cup} R_m$.
\end{lemma}

\noindent
{\bf Proof:} Since $R_0^+ \stackrel{\cdot}{\cup} R_m \subset
R^+$ it is sufficient to check that the cardinalities coincide,
i.e.\ $|R| - |R_0| = 2|R_m|$. We have the natural isomorphisms
\[ T_pM \cong {\fr g}^0/{\fr h}^0 \cong {\fr g}/{\fr g}_p\, ,
\quad (T_pM)\otimes \mbox{\C} \cong {\fr g}/{\fr h}\]
and hence
\[ |R| - |R_0| = \dim_{\mbox{\C}}{\fr g}/{\fr h}
= \dim_{\mbox{\R}}{\fr g}/{\fr g}_p = 2(|R| - |R_0|) -2|R_m|\, .  \]
This proves Lemma \ref{parabLemma} and shows that ${\fr g}_p = {\fr a} +
\sum_{\alpha \in R_0\cup R_m} {\fr g}_{\alpha} \subset {\fr g}$ is a parabolic
subalgebra. Let $P \subset G$ be the parabolic subgroup with Lie algebra
${\fr p} = {\fr g}_p$. The intersection $P\cap G^0$ is a Lie group
with Lie algebra ${\fr p} \cap {\fr g}^0 = {\fr h}^0$. Therefore,
$H^0 \subset P\cap G^0$ is the identity component of $P\cap G^0$  and
\[ P\cap G^0 = N_P(H^0)\cap G^0 = N_P(Z_G(h))\cap G^0 = Z_G(h)\cap G^0 = H^0
\, .\]
We have used that $P = Z_G(h) \mbox{\Bbb n} N$ is the semidirect product
of the maximally reductive subgroup $Z_G(h)$ with the nilpotent radical
$N$ and hence $Z_G(h)$ equals its own normalizer in $P$.
This completes the proof of Proposition \ref{basicProp}. 

Now we specialize Theorem \ref{redThmI}:
\begin{thm} \label{redThmII} Let $(M,g)\subset P(V)$ be a
homogeneous special pseudo-K\"ahler
manifold
of the connected real semisimple group $G^0$ with compact stabilizer
$H^0 = G^0_p$, $p\in M$. Then the stabilizer $H^0 = Z_{G^0} (h)$ is the
centralizer of some element $h \in {\fr h}^0 = Lie \, H^0$ and is connected;
$M = G^0/H^0$ is simply connected. The orbit $\overline{M} = Gp = GM$ under
the complexified linear group $G = (G^0)^{\mbox{\C}} \subset GL(V)$ is
a flag manifold and the stabilizer $P = G_p \subset G$ is a parabolic
subgroup with maximally reductive subgroup $H = (H^0)^{\mbox{\C}}$.
Finally, $\overline{M}$ is precisely the Zariski closure of $M$
in $P(V)$. In particular, we have the canonical open inclusion
$M \subset \overline{M}$.
\end{thm}

\noindent
{\bf Proof:} We denote by $J$ the $G^0$-invariant complex structure on $M$
induced by the inclusion $M \subset P(V)$.
The first part of Theorem \ref{redThmII} follows immediately from 
Theorem \ref{redThmI}. Moreover, since $G^0$ acts effectively on $M$, the
Lie group $G^0$ is isomorphic to its adjoint group $Ad \, G^0$.
We will write $G = (G^0)^{\mbox{\C}}$. Now Proposition \ref{basicProp}
shows that the stabilizer $P = G_p$ is a parabolic subgroup
(with maximally reductive subgroup $H = (H^0)^{\mbox{\C}}$) and therefore
$Gp \cong G/P$ a flag manifold. In particular, $\overline{M} = Gp
\subset P(V)$ is a compact complex submanifold and hence Zariski closed
(by Chow's Theorem). $\overline{M}$ being connected and $\dim M = \dim
\overline{M}$ this implies that $M \subset \overline{M}$ is open
(in the standard Hausdorff topology) and $\overline{M}$ its Zariski
closure. This proves Theorem \ref{redThmII}.  

\section{Classification of Lagrangian cones over flag ma\-ni\-folds}
If $M \subset P(V)$ is a homogeneous special pseudo-K\"ahler manifold
as in Theorem \ref{redThmII} then, by virtue of this theorem, its Zariski
closure $\overline{M}\subset P(V)$ is an orbit of the complex semisimple
group $G \subset Sp(V)$. The stabilizer of any point $p\in \overline{M}$
is a parabolic subgroup of $G$ and hence $p$ is a {\bf highest weight
direction} of the $G$-module $V$, i.e.\ $p = \mbox{\C}^{\ast}v$ for some
highest weight vector $v\in V$ (with respect to some choice of Cartan
subalgebra ${\fr a} \subset {\fr  g}$ and system of positive roots $R^+$).
The cone ${\cal C}_{\overline{M}} = \{ v\in V| \pi v \in \overline{M}\}$
($\pi : V \rightarrow P(V)$ the canonical projection) over the
Zariski closure  $\overline{M}$
is Lagrangian, since the cone ${\cal C}_M$ over $M$ is Lagrangian.
Moreover, $G$ acts
transitively on ${\cal C}_{\overline{M}}$. In fact, ${\cal C}_{\overline{M}}$
consists of highest weight vectors and hence for all
$v\in {\cal C}_{\overline{M}}$ we can find an algebraic torus in $G$ acting
transitively on the direction $\mbox{\C}^{\ast} v$ spanned by $v$.
The claim follows now from the fact that $G$ acts transitively on the set
of directions $\overline{M}$. In this section we will classify all
symplectic modules
$V$ of complex semisimple groups $G$ which contain a highest weight vector
$v$ whose orbit $Gv\subset V$ is Lagrangian. Propositions
\ref{irredProp} and \ref{simpleProp} reduce this problem to the case of
irreducible modules of simple groups.

\begin{prop} \label{irredProp}  Let $(V,\omega )$ be a symplectic $G$-module
of a connected complex semisimple group $G$ and ${\cal C} = Gv$ a Lagrangian
$G$-orbit of a highest weight vector $v\in V$. Then either
\begin{enumerate}
\item[i)] $V = U \oplus U^{\ast}$ is the direct sum of an irreducible
$G$-module $U$ and its dual $U^{\ast}$ and we have $U = {\cal C} \cup
\{ 0\}$ or
\item[ii)] $V$ is irreducible.
\end{enumerate}
In the first case ${\fr g} = Lie \, G = {\fr sl}_n(\mbox{\C})$ or
${\fr sp}_n(\mbox{\C})$
and $U$ is the standard $\fr g$-module, i.e.\ $U = \mbox{\C}^n$ 
if ${\fr g} = {\fr sl}_n(\mbox{\C})$  and $U = \mbox{\C}^{2n}$ 
if ${\fr g} = {\fr sp}_n(\mbox{\C})$.
\end{prop}

\noindent
{\bf Proof:} The cone $\cal C$ is contained in a unique irreducible submodule
$U\subset V$. By Schur's Lemma the restriction $\omega |U$ of the symplectic
form $\omega$ to the irreducible submodule $U \subset V$ is either zero or
nondegenerate. In the first case $U$ is totally isotropic and contains the
Lagrangian cone $\cal C$ as an open orbit. On the other hand, it is 
known, see \cite{26}, that  the orbit $Gp$ of the highest weight 
direction $p =  \mbox{\C}^{\ast}v$ is closed in $P(V)$. 
Thus $U = {\cal C} \cup \{ 0\}$ is 
a Lagrangian subspace .
Let $U'$ be a $G$-invariant complement to $U$ in $V$. The symplectic structure
$\omega$ defines an isomorphism $U' \cong U^{\ast}$ of $G$-modules. {}From
the fact that $G$ is a complex semisimple Lie group acting transitively
on the projective space $P(U)$ by projective transformations it follows that
${\fr g} = {\fr sl}_n (\mbox{\C})$ or ${\fr g} = {\fr sp}_n (\mbox{\C})$
and $U$ is the standard ${\fr g}$-module (up to an automorphism of $\fr g$),
see \cite{20}. This completes the discussion of the first case. In the
second case  ${\cal C}$ is Lagrangian in $(U,\omega |U)$
and hence $V = U$. In particular, $V$ is irreducible. This proves
the proposition. 

Recall that any connected and simply connected complex semisimple Lie
group $G$ is the direct product $G = G_1 \times \cdots \times G_k$
of simple groups $G_i$. Any irreducible $G$-module $V$ is the tensor
product $V = V_1 \otimes \cdots \otimes V_k$ of irreducible $G_i$-modules
$V_i$.

\begin{lemma} \label{symplecticLemma} The irreducible $G$-module
$V = V_1 \otimes \cdots \otimes V_k$
is symplectic if and only if each of the $G_i$-modules $V_i$ is self-dual
and the number of symplectic factors $V_i$ is odd.
\end{lemma}

\noindent
{\bf Proof:} We choose a Cartan subalgebra ${\fr a}_i \subset {\fr g}_i
= Lie \, G_i$ and a system of positive ${\fr a}_i$-roots $R^+_i$ for each of
the simple
summands ${\fr g}_i$ of ${\fr g} = Lie \, G$. Let $\Lambda^i$
(respectively $M^i$) be
the highest (respectively lowest) weight of the $G_i$-module $V_i$ with
respect to $R^+_i$. Then
${\fr a} = \oplus_{i = 1}^k {\fr a}_i \subset {\fr g}$ is a Cartan
subalgebra and $R^+ = \cup_iR^+_i$ a system of positive ${\fr a}$-roots
for $\fr g$, where we
have identified ${\fr a}_i^{\ast} = \{ \alpha \in {\fr a}^{\ast}|
\alpha (a) = 0$ for all $a\in {\fr a}_j$, $j\neq i\}$.
The $G$-module $V$ has highest weight $\Lambda = \sum_{i=1}^k \Lambda^i$ and
lowest  weight $M = \sum_{i=1}^k M^i$ with respect to $R^+$.
$V$ is self-dual if and only if $\Lambda + M = 0$, see \cite{21}.
Restricting $\Lambda + M$ to ${\fr a}_i$ we see that $V$ is self-dual if
and only if $\Lambda^i + M^i = 0$ for all $V_i$, i.e.\ if and only if
all $V_i$ are self-dual. Let $b_i\neq 0$ be a $G_i$-invariant symmetric or
skew symmetric bilinear form on $V_i$.
Then $b = b_1 \otimes \cdots \otimes b_k$ is a $G$-invariant bilinear
form on $V$. It is skew symmetric if and only if the number of
skew symmetric factors $b_i$ is odd. This proves the lemma. 

The next proposition shows, in particular, that the highest weight vector
orbit $\cal C$ of an irreducible symplectic module $V$ of a simple group
is automatically Lagrangian provided that $\dim V = 2\dim {\cal C}$.

\begin{prop}\label{totisotrProp} Let $(V,\omega )$ be an irreducible
symplectic module
of a complex semisimple group $G\subset Sp(V)$ and ${\cal C} \subset V$ its
highest weight vector orbit. Then $\cal C$ is a totally isotropic cone,
unless $V = V_1 \otimes \cdots \otimes V_k$ is the tensor product of an odd
number $k$ of symplectic vector spaces $(V_i,\omega_i)$, $\omega =
\omega_1 \otimes \cdots \otimes \omega_k$,  and $G = Sp(V_1) \otimes \cdots \otimes
Sp(V_k)$ is the image of the group $Sp(V_1)\times \cdots \times Sp(V_k)$
under the tensor product representation on $V$. If $(V,g)$ is an orthogonal
module of a complex semisimple group $G$ then any highest weight vector
orbit is contained in the null cone of $(V,g)$.
\end{prop}

\noindent
{\bf Proof:} The stabilizer of a point $p\in P({\cal C})$ being parabolic,
there exists a Cartan subalgebra ${\fr a} \subset {\fr g} = Lie \, G$ and a
system of positive roots $R^+$ (of $\fr g$ with respect to $\fr a$) such
that $\sum_{\alpha \in R^+} {\fr g}_{\alpha} v = 0$ and $av = \Lambda (a) v$
for all $a \in {\fr a}$, where $\Lambda$ is the highest weight of the
irreducible $\fr g$-module $V$. In particular, we have
\[ T_v{\cal C} = {\fr g} v = \mbox{\C}v + \sum_{\alpha \in R^-}
{\fr g}_{\alpha}v\, ,\quad R^- = - R^+\, .\]
It is sufficient to show that the tangent space $T_v{\cal C}$ is totally
isotropic, unless ${\fr g} = sp(V_1) \oplus \cdots \oplus sp(V_k)$. For
$a\in {\fr a}$, ${\fr g}_{\alpha} = \mbox{\C}e_{\alpha}$, $\alpha \in R^-$,
using the $G$-invariance of $\omega$ we compute
\[ \Lambda (a) \omega (v,e_{\alpha}v) = \omega (av,e_{\alpha}v)
= -\omega (v,[ a,e_{\alpha}] v + e_{\alpha}av) =
-(\alpha (a) + \Lambda (a)) \omega (v,e_{\alpha}v)\]
and hence $\omega (v,e_{\alpha}v) = 0$, unless $-\alpha = 2\Lambda$.

Assume for the moment that $-\alpha \neq 2\Lambda$. By our computation,
this implies that $\omega (v,{\fr g}v) = 0$.
In order to prove that $T_v {\cal C} = {\fr g} v$ is totally isotropic
it only remains to show that $\omega (e_{\alpha}v,e_{\beta}v) = 0$ for
all $\alpha ,\beta \in R^-$. Using the fact that $\omega (v ,{\fr g} v) = 0$,
we compute
\[ \omega (e_{\alpha}v,e_{\beta}v) = - \omega (v,e_{\alpha}e_{\beta}v) =
- \omega (v,[e_{\alpha},e_{\beta}] v) - \omega (v,e_{\beta}e_{\alpha}v) \]
\[ = 0 + \omega (e_{\beta}v,e_{\alpha}v) = - \omega (e_{\alpha}v,e_{\beta}v)
\, ,\]
which proves that $\omega (e_{\alpha}v,e_{\beta}v) = 0$.

To discuss the case $-\alpha = 2\Lambda$ let ${\fr g} = {\fr g}_1
\oplus \cdots \oplus {\fr g}_k$ be the decomposition of $\fr g$ into simple
ideals and ${\fr a} = {\fr a}_i \oplus \cdots \oplus {\fr a}_k$ the
corresponding decomposition of the Cartan subalgebra
${\fr a} \subset {\fr g}$ into Cartan subalgebras
${\fr a}_i \subset {\fr g}_i$. Since $G$ acts faithfully on $V$ we have
$\Lambda |{\fr a}_i \neq 0$. If $-\alpha = 2\Lambda$ then $-\alpha |{\fr a}_i
\neq 0$ is a dominant root of the simple Lie algebra ${\fr g}_i$,
i.e.\ $\langle -\alpha |{\fr a}_i,\beta |{\fr a}_i\rangle \ge 0$ for all
$\beta \in R^+$. This implies that $-\alpha |{\fr a}_i = \delta ({\fr g}_i)$
is the highest root of ${\fr g}_i$. Moreover, $\delta ({\fr g}_i)/2 =
\Lambda |{\fr a}_i$ is an element of the weight lattice, i.e.\
is an integral linear combination of fundamental weights. This does only
occur if ${\fr g}_i$ is the symplectic Lie algebra and $\Lambda |{\fr a}_i$
is the first fundamental weight, which is the highest weight of the
standard representation. This proves the proposition for the symplectic case.
In the orthogonal case it is sufficient to show that the highest weight
vector $v$ is a null vector. This follows from the computation:
\[ \Lambda (a) g(v,v) = g(av,v) = - g( v,av) = -\Lambda (a) g(v,v)\, ,
\quad a\in {\fr a} \, .  \]
This finishes the proof of Proposition \ref{totisotrProp}. 

\begin{prop} \label{simpleProp} Let $(V,\omega )$ be an irreducible
symplectic module of a connected complex semisimple group $G\subset Sp(V)$
and ${\cal C} \subset V$ its highest weight vector orbit. Assume that
$\cal C$ is Lagrangian. Then either $G$ is simple or
$V = \mbox{\C}^2\otimes \mbox{\C}^n$ and   $G = SL_2(\mbox{\C}) \otimes
G'$, where $\mbox{\C}^2$ is the standard module of $SL_2(\mbox{\C})$ and
$\mbox{\C}^n$ is the standard module of $SO_n(\mbox{\C})$ considered
as module of a simple subgroup $G'\subset SO_n(\mbox{\C})$ acting
transitively on the cone ${\cal C}' \subset \mbox{\C}^n$ of nonzero
null vectors. More precisely, $G' = SO_n(\mbox{\C})$ or
$G' = G_2(\mbox{\C})\subset SO_7(\mbox{\C})$. Moreover, $P({\cal C}) =
P(\mbox{\C}^2) \times Q
\subset P(\mbox{\C}^2\otimes \mbox{\C}^n) = P(V)$ is
the product of the projective line $P(\mbox{\C}^2)$ with the quadric
$Q = P({\cal C}') \subset P(\mbox{\C}^n)$.
\end{prop}

\noindent
{\bf Proof:} If $G$ is not simple, by Lemma \ref{symplecticLemma}, we can
write $V = V_1 \otimes V'$ where $(V_1,\omega_1)$
is an irreducible symplectic module of a simple group $G_1$,
$(V',g')$ an irreducible orthogonal module of a semisimple group $G'$
and $G = G_1 \otimes G'$.  Let $v \in {\cal C}$, then $v = v_1\otimes v'$,
where $v_1$ (respectively $v'$) is a highest weight vector of the
$G_1$-module $V_1$ (respectively of the $G'$-module $V'$). The tangent
space
\[ T_v{\cal C} = {\fr g} v = {\fr g}_1v_1 \otimes v' + v_1 \otimes
{\fr g}'v'\]
(${\fr g} = Lie \, G$, ${\fr g}_1 = Lie \, G_1$, ${\fr g}' = Lie \, G'$) is
Lagrangian. By Proposition 
\ref{totisotrProp}, $v' \in V'$ is a null vector for the
metric $g'$ and hence ${\fr g}'v'$ is tangent to the null cone ${\cal C}'
\ni v'$. This implies that $g(v',{\fr g}'v') = 0$. Now a straightforward
calculation shows that $V_1 \otimes v' + v_1 \otimes T_{v'}{\cal C}'$ is
totally isotropic for the symplectic form $\omega = \omega_1 \otimes g'$.
The subspace $T_v{\cal C} \subset V_1\otimes v' + v_1 \otimes T_{v'}{\cal C}'$
being Lagrangian, it follows that ${\fr g}_1v_1 = V_1$, ${\fr g}'v' =
T_{v'}{\cal C}'$ and $2(\dim V_1 +
\dim V' -2) = \dim V_1 \cdot \dim V'$. Since $G$ is semisimple, 
this is only possible if
$\dim V_1 = 2$ and hence $G_1 = Sp_1(\mbox{\C}) = SL_2(\mbox{\C})$.
The remaining statements follow from the classification of transitive
projective actions of complex semisimple groups on the quadric
$Q = P({\cal C}') \subset P(\mbox{\C}^n)$, see \cite{20}. 
This establishes the proposition.

Now we associate with a complex simple Lie algebra $\fr l$ a
symplectic $G$-module $V$ such that the orbit of the highest
weight vector is a Lagrangian cone.

  Let $\mu$ be the highest root of $\fr l$ (with respect to
some Cartan subalgebra  and choice of positive root system).
Denote by
$${\fr g}(\mu) = {\rm span}\{ h_{\mu}=[e_{\mu}, \, e_{-\mu}],\,
e_{\mu},\, e_{-\mu}\} \cong {\fr sl}_2(\mbox{\Bbb C}) $$
the corresponding 3-dimensional subalgebra, where the root vectors
$e_{\pm \mu}$ are nor\-ma\-li\-zed by the condition
$$ \langle e_{\mu}, \, e_{- \mu} \rangle = \frac{2}{\langle \mu,\, 
\mu \rangle }, $$
 (and $\langle \cdot , \cdot \rangle$  is the Killing form).
 One can check easily that the adjoint operator ${\rm
ad}(h_{\mu})$ has eigenvalues $ 0, \pm 1, \pm 2 $ and the
corresponding eigenspace decomposition
$$ {\fr l} = {\fr l}_{-2} +{\fr l}_{-1} + {\fr l}_0+ {\fr l}_1 
+ {\fr l}_2
$$
defines a grading on the Lie algebra $\fr l$. Moreover,
${\fr l}_{\pm 2} = \mbox{\Bbb C}e_{\pm \mu} $ and
${\fr l}_0 = \mbox{\Bbb C}h_{\mu} + {\fr l}'_0 $
where ${\fr l}'_0$ is the centralizer of ${\fr g}(\mu)$.

\noindent 
  Denote by $  L^{\tau} $ the compact Lie group without centre
associated with $\fr l$ and by
   $$L^{'\tau}_0, \qquad L^{\tau}_0= T^1\cdot L^{'\tau}_0,
\qquad L_{ev}^{\tau}  = Sp_1 \cdot L^{'\tau}_0 $$
its compact subgroups  associated with the  subalgebras
 $$  {\fr l}'_0,\qquad {\fr l}_0\quad \mbox{and} \quad
  {\fr l}_{ev} =  {\fr l}_{-2} + {\fr l}_0
+ {\fr l}_2 \quad \mbox{respectively}.$$
 The following result is known.

\begin{thm}
\begin{enumerate}
\item $M= L^{\tau}/L_{ev}^{\tau}$ is a compact symmetric quaternionic
K\"ahler ma\-ni\-fold (called Wolf space). Moreover, the Wolf
spaces exhaust all homogeneous quaternionic K\"ahler
manifolds with positive scalar curvature.
\item
 $$ Z = L^{\tau}/L_0^{\tau} \longrightarrow M = L^{\tau}
/L_{ev}^{\tau} $$
is the twistor fibration of the quaternionic K\"ahler manifold
$M$. In particular, the manifold $Z$ has the structure of a
complex contact homogeneous manifold and admits an invariant 
K\"ahler-Einstein metric, which is unique up to scaling. Moreover, such manifolds $Z$
exhaust all compact homogeneous  complex contact manifolds.
\end{enumerate}
\end{thm}

  Denote by ${\fr g} = [ {\fr l}_0 , {\fr l}_0] \subset {\fr l}_0'$ the 
Levi subalgebra of ${\fr l}_0$ and by $G$ the corresponding 
complex Lie subgroup of $Aut ({\fr l})$. Note that ${\fr g} = 
{\fr l}_0'$, unless ${\fr l} = {\fr sl}_{n+2}(\mbox{\C})$; 
${\fr g} = {\fr sl}_n (\mbox{\C}) \subset {\fr l}_0' = 
{\fr gl}_n (\mbox{\C})$ if ${\fr l} = {\fr sl}_{n+2} (\mbox{\C})$. 
The adjoint representation of $G$ in the
space $V = {\fr l}_{-1} $ is symplectic, with the symplectic form
defined by the Lie bracket $[\cdot , \cdot ]: \wedge^2{\fr l}_{-1}
\rightarrow {\fr l}_{-2}$.
We will call $V$ the {\bf standard symplectic $G$-module
associated with the simple Lie algebra} $\fr l$.

In the following table we indicate for any complex simple Lie
algebra $\fr l$ the associated Wolf space $N =
L^{\tau}/L_{ev}^{\tau} $ , complex Lie group $G$ ,
symplectic $G$-module $V = {\fr l}_{-1} $ and  the
semisimple part $H'$ of the stabilizer of a highest weight vector
of $V$. 
We denote by $(V^{(d)} = V(\Lambda ), {\fr g})$   the
irreducible $\fr g$-module with highest weight $\Lambda = \sum
\Lambda_i \pi_i$, $\pi_i$ the fundamental weights, and the
superindex $d$ indicates the dimension.

Table 1\\ 
\begin{tabular}{|l|l|l|l|l|}
\hline
 ${\fr l}$ & $N$ & $G$ & $V$ & $H'$ \\
\hline
$ {\fr sl}_{n+2}(\mbox{\Bbb C})$ &
$ SU_{n+2}/S(U_2\cdot U_n)$ &
$ SL_n(\mbox{\Bbb C})$ &
$ V^{(2n)}=\mbox{\Bbb C}^n \oplus \mbox{\Bbb C}^{n*}$ &
$ SL_{n-1}(\mbox{\Bbb C})$ \\
\hline
$ {\fr so}_{n+4}(\mbox{\Bbb C})$ &
$ SO_{n+4}/SO_4 \cdot SO_n $ &
$ SL_2(\mbox{\C}) \cdot SO_n(\mbox{\C})$ &
$ V^{(2n)} =\mbox{\C}^2 \otimes \mbox{\C}^n $ &
$ SO_{n-2}(\mbox{\C})$ \\
\hline
${\fr sp}_{n+1}(\mbox{\C})$ &
$ Sp_{n+1}/Sp_1\cdot Sp_n $&
$Sp_{n}(\mbox{\C})$ &
$ V^{(2n)} =\mbox{\C}^{2n}$ &
$ Sp_{n-1}(\mbox{\C})$ \\
\hline
${\fr g}_2(\mbox{\C})$ & $ G_2/SU_2\cdot SU_2 $ & $SL_2(\mbox{\C})$ &
$V^{(4)} =V(3\pi_1)$ & $\{ e\}$ \\
\hline
${\fr f}_4(\mbox{\C}) $ & $F_4/Sp_1 \cdot Sp_3$ &$ Sp_3(\mbox{\C})$ &
 $V^{(14)} =V(\pi_3)$ &$ SL_3(\mbox{\C})$ \\
\hline
${\fr e}_6(\mbox{\C}) $ & $ E_6/Sp_1 \cdot SU_6$ & $SL_6(\mbox{\C})$ &
$V^{(20)} = V(\pi_3) =\wedge^3\mbox{\C}^6 $ &
 $SL_3(\mbox{\C})\cdot SL_3(\mbox{\C})$ \\
\hline
${\fr e}_7(\mbox{\C}) $ & $ E_7/Sp_1 \cdot Spin_{12} $ &
 $ Spin_{12}(\mbox{\C})$ & $ V^{(32)} = V(\pi_5) $ &
$SL_6(\mbox{\C})$ \\
\hline
$ {\fr e}_8(\mbox{\C}) $ & $ E_8/Sp_1 \cdot E_7 $ & $E_7(\mbox{\C})$ &
$ V^{(56)} = V(\pi_1)$ & $E_6(\mbox{\C}) $ \\
\hline
\end{tabular}

 The following theorem shows that with a few number of exceptions
 the standard symplectic modules exhaust all $G$-modules with
 Lagrangian orbit of the highest weight vector $v$.

\begin{thm} \label{listthm}
\begin{enumerate}
\item Let $V$ be the standard $G$-module associated with a simple
Lie algebra ${\fr l} \neq \fr sp_n (\mbox{\C}) $.  Then the orbit
${\cal C} = Gv$ of the highest weight vector $v$ is a Lagrangian cone.
\item  Let $V$ be a symplectic $G$-module of a semisimple Lie group $G$
such that the orbit of a highest weight vector $v$ is a Lagrangian cone.
Assume that $V$ is not a standard symplectic module. Then the structure
of $G$-module on $V$ can be extended to the structure of a standard symplectic
$\tilde G$-module $V$, where $\tilde G$ is a semisimple Lie group which
contains $G$ , such that the $G$-orbit and the $\tilde G$-orbit of the
highest weight vector $v$ coincide. More precisely, this happens 
only in the following three cases:
\end{enumerate} 
\begin{enumerate} 
\item[i)] $G= Sp_n(\mbox{\C}) \subset \tilde G = SL_{2n}(\mbox{\C}), 
\qquad V = \mbox{\C}^{2n}  \oplus \mbox{\C}^{2n}{}^*, \qquad
   H' = Sp_{n-1}(\mbox{\C}) $

\item[ii)] $G = Spin_{11}(\mbox{\C}) \subset \tilde G = Spin_{12}(\mbox{\C}), 
\qquad  V = V^{(32)} = V(\pi_5) $  (the semispinor module of $Spin_{12}
(\mbox{\C}) $, which is the spinor module of $ Spin_{11}(\mbox{\C})),
 \qquad H' = SL_5 (\mbox{\C}) $

\item[iii)] $G= SL_2 (\mbox{\C}) \times G_2(\mbox{\C}) \subset \tilde G =
SL_2(\mbox{\C}) \times SO_7(\mbox{\C}), \qquad V = V^{(14)} = V(\pi_1 ) \otimes V(\pi_1 ) =\mbox{\C}^2 \otimes
 \mbox{\C}^7, \qquad H' = SL_2 (\mbox{\C}) . $
\end{enumerate}
(Here $H'$ is the semisimple part of the stabilizer of $v$ in $G$).

\end{thm}
\noindent
{\bf Proof:} First we prove 1). In the case ${\fr l} = {\fr sl}_{n+2}
(\mbox{\C})$ the highest weight vector orbits $\mbox{\C}^n -\{ 0\}$
and  $(\mbox{\C}^n)^{\ast} -\{ 0\} \subset \mbox{\C}^n \oplus
(\mbox{\C}^n)^{\ast} = V$ are Lagrangian orbits of the group
$G = SL_n(\mbox{\C})$. In the remaining cases $V$is irreducible
and by  Proposition \ref{totisotrProp} and Lemma  \ref{dimLemma} the $G$-orbit
$\cal C$ of the highest weight vector $v$ is Lagrangian if and only if 
$$ 2 \dim {\cal C } = \dim V .$$
  The dimension of $\cal C$ can be calculated by the formula
  \[ \dim {\cal C} = \frac{1}{2} (\dim G - \mbox{rk}\,  G - \dim H' +
\mbox{rk}\,  H') + 1\, ,\]
where  $H'$ is the  semisimple part of the stabilizer of $v$.
Indeed,
 $\dim P({\cal C})  = |R^+| - |R^+_0|$,
  where $R^+$ (respectively $R^+_0$) is a system of
positive roots for $\fr g$ (respectively ${\fr h}'$). Since $|R^+| =
(\dim G - \mbox{rk}\,  G)/2$ and $|R^+_0| =
(\dim H'- \mbox{rk}\,  H')/2$ we  get this formula.
 Using it  and the table, one can check immediately that for each
 standard symplectic $G$-module the orbit of the highest weight vector
 is Lagrangian.

 Now we prove 2). By  Proposition \ref{irredProp} and Proposition 
\ref{simpleProp} 
it is sufficient to consider the case when  $V$ is an irreducible 
symplectic module of a simple Lie group $G$.
 The proof is based on the following lemma and proposition.

 \begin{lemma} \label{dimLemma} Let $V$ be a symplectic module of a complex
semisimple
group $G$ such that the highest weight vector orbit ${\cal C} \subset V$
is Lagrangian. Then
\[ \dim V \le \dim G  - \mbox{{\rm rk}}\, G + 2\, .\]
\end{lemma}

\noindent
{\bf Proof:} Since $P({\cal C}) \subset P(V)$ is a flag manifold of the
group $G$, we have $\dim P({\cal C}) \le (\dim G - \mbox{rk}\, G )/2$ and
hence
\[ \dim V = 2 \dim {\cal C} = 2(\dim P({\cal C}) + 1) \le
\dim G  - \mbox{rk}\, G + 2\, , \]
proving the lemma. 

 \begin{prop} The irreducible symplectic $G$-modules of  complex simple
 (connected) Lie group $G$ which satisfy the condition
 $$\dim V \leq \dim G - {\rm rk}G +2 $$
 are listed below. If the module is standard symplectic, we indicate also
 the corresponding simple Lie algebra $\fr l$.
\begin{enumerate}
\item[A)] $ G = SL_2(\mbox{\C}), \quad   V^{(4)} = V(3\pi_1) = \vee^3 
\mbox{\C}^2,\quad     {\fr l}= {\fr g}_2 (\mbox{\C})  $,\\
$ G = SL_6 (\mbox{\C}) , \quad V^{(20)} = V(\pi_3) = \wedge^3  
\mbox{\C}^6, \quad {\fr l} = {\fr e}_6 (\mbox{\C}) $,\\
\item[B)] $ G =Spin_{11}(\mbox{\C}), \quad
  V^{(32)} = V(\pi_5) =$ spinor module, \\
$G =Spin_{13}(\mbox{\C}), \quad
V^{(64)} = V(\pi_6) =$ spinor module,\\
\item[C)] $G = Sp_l(\mbox{\C}), \quad V^{(2l)} = V(\pi_1)
 = \mbox{\C}^{2l},\quad {\fr l}={\fr sp}_{l+1}(\mbox{\C}),$ \\
$G = Sp_3(\mbox{\C}), \quad V^{(14)} =  V(\pi_3), \quad
{\fr l} ={\fr f}_4(\mbox{\C})$,\\

\item[D)] $G = Spin_{12}(\mbox{\C}), \quad
V^{(32)} = V(\pi_5$ or $\pi_6) =$ semi spinor modules,
$ \quad {\fr l } ={\fr e}_7(\mbox{\C})$,
\item[E)] $ G = E_7(\mbox{\C}), \quad  V^{(56)} = V(\pi_1),\quad
 {\fr l} = {\fr e}_8(\mbox{\C})$.
\end{enumerate}
\end{prop}

 \noindent
{\bf Proof:} The proof is a straightforward (but long) exercise; the
needed representation theory and useful tables can be found in
\cite{21} and \cite{22}. 
\noindent

Now to finish the proof of the theorem, it only remains to investigate 
the two modules in the case B, which are the only nonstandard 
symplectic modules from the list.

For the spinor module $V^{(64)}$ of  $Spin_{13}(\mbox{\C})$ 
  the semisimple part of the stabilizer of the highest weight vector 
  $v$ is $SL_6(\mbox{\C})$ . 
This can be easily derived from the fact that the Dynkin diagram
 of the root system $R_0$ of ${\fr h}' = {\rm Lie} \, H'$ is obtained
from the Dynkin diagram  of  the Lie algebra $\fr g$ of $G$ 
by deleting all simple roots
$\alpha_i$ such that $\langle \alpha_i, \Lambda \rangle \neq 0$.
Hence,
\[ \dim {\cal C} = \frac{1}{2} (\dim G - \mbox{rk}\,  G - \dim H' +
\mbox{rk}\,  H') + 1\ = \]
\[ =  \frac{1}{2} (6\cdot 13 - 6 -35 + 5) + 1 = 22 <
\frac{1}{2}\dim V = 32 \]
and the orbit $\cal C $ is not Lagrangian in this case. 

The restriction of the semispinor module of
$Spin_{12}(\mbox{\C})$ to a module of $Spin_{11}(\mbox{\C}) \subset
Spin_{12}(\mbox{\C})$ gives the spinor module of $Spin_{11}(\mbox{\C})$.
This implies that the highest weight vector orbit of the group
$Spin_{11}(\mbox{\C})$ is contained in the highest weight vector orbit of
$Spin_{12}(\mbox{\C})$. A short dimension count as above
shows that the two orbits coincide and, hence, are 
Lagrangian cones. This establishes the proposition.

\noindent
{\bf Remark 2:} For any of the complex contact manifolds $Z$ listed in 
Theorem 3 and Table 1 with the exception of $Z= 
SU_{n+2}/U_1 \cdot U_n = P(T^*\mbox{\C}P^{n+1})$ the
second Betti number is 1. Hence, using the exact sequence
$0 \rightarrow \mbox{\Z} \rightarrow {\cal O} \rightarrow {\cal O}^{\ast}
\rightarrow 0$, one can easily show that the Picard group $Pic(Z) =
H^1 (Z,{\cal O}^{\ast}) \cong \mbox{\Z}$. Let $\cal L$ be the ample
generator of $Pic(Z)$. The {\bf degree} $k$ of a rational
curve $\mbox{\C}P^1 \hookrightarrow Z$ is defined by
\[  {\cal L}|\mbox{\C} P^1 \; \cong \; {\cal O} (k) \]
and we can consider the cone ${\cal C}_{min} \subset  T_oZ$ of nonzero tangent
vectors to rational curves of minimal degree $k = 1$ through the point
$o \in Z$. If $Z =\mbox{\C} P^{2n+1}$ then ${\cal C}_{min} = T_oZ
- \{ 0\}$. In the remaining cases ${\cal C}_{min}$ is contained in the 
contact hyperplane $W \subset  T_oZ$ and is precisely the highest weight
vector orbit of the $G$-module $W$, where $G$  is the complexification
of the Levi subgroup of the isotropy group  $L_0^{\tau}$ of 
$Z = L^{\tau}/L_0^{\tau}$, cf.\ \cite{16}.

\section{Real forms} By Theorem \ref{redThmII} any homogeneous special 
pseudo-K\"ahler manifold of a real semisimple group $G^0$ 
with compact stabilizer
is an open orbit of $G^0$ on a flag manifold $P({\cal C})$, where
${\cal C} \subset V$ is a Lagrangian highest weight vector orbit
in a symplectic $G$-module $(V,\omega )$, $Lie \, G =
(Lie \, G^0) \otimes \mbox{\C}$. Moreover, $V$ carries the $G^0$-invariant
real structure $\tau$ and Hermitian form $\gamma = \sqrt{-1} \omega (\cdot ,
\tau \cdot )$. In Theorem \ref{listthm} we have 
enumerated all symplectic modules $V$ of
complex semisimple group $G$ for which a highest weight vector orbit is
Lagrangian. Now we list all real forms of those modules.

\begin{prop} \label{realformProp} For any symplectic module
$(V,\omega )$ of a complex
semisimple Lie algebra $\fr g$ with Lagrangian highest weight vector orbit
the list of real forms ${\fr g}^0$ of ${\fr g}$ for which there
exists a ${\fr g}^0$-invariant real structure on $V$ is given below (up to
an automorphism of $\fr g$). For the simple Lie algebras $\fr g$ we
have:
\begin{enumerate}
\item[A)] $\vee^3 \mbox{\C}^2$:  ${\fr sl}_2 (\mbox{\R})$.\\
$\wedge^3  \mbox{\C}^6$:  ${\fr sl}_6 (\mbox{\R})$, ${\fr su}_{1,5}$ and
$ \quad {\fr su}_{3,3}$.\\
$\mbox{\C}^n \oplus (\mbox{\C}^n)^{\ast}$: ${\fr sl}_n (\mbox{\R})$ and
$ \quad {\fr su}_{k,l}$, \quad $k+l = n$.
\item[B)] Spinor module: ${\fr so}_{1,10}$ and $ \quad {\fr so}_{2,9}$.
\item[C)] $V(\pi_3)$: ${\fr sp}_3(\mbox{\R})$.\\
$\mbox{\C}^{2n}\oplus (\mbox{\C}^{2n})^{\ast} \cong 2V(\pi_1)$:
${\fr sp}_n(\mbox{\R})$ and  ${\fr sp}(k,l)$, $k+l = n$.
\item[D)] $V(\pi_5)$: ${\fr so}_{2,10}$ and $ \quad {\fr so}_{6,6}$.\\
$V(\pi_6)$: ${\fr so}_{2,10}$, ${\fr so}_{6,6}$ and ${\fr so}_{12}^{\ast}$.
\item[E)] $V(\pi_1)$: ${\fr e}_7^{(i)}$, $i = -25,-5$ or $7$.
\end{enumerate}
The remaining cases are:
$\mbox{\C}^2 \otimes \mbox{\C}^n$: ${\fr sl}_2(\mbox{\R}) \oplus
{\fr so}_{k,l}$ and, if $n$ is even, also ${\fr su}_2 \oplus
{\fr so}_n^{\ast}$.\\
$\mbox{\C}^2 \otimes \mbox{\C}^7$: ${\fr sl}_2(\mbox{\R}) \oplus
{\fr g}_2^{(2)}$.
\end{prop}

\noindent
{\bf Proof:} The proof is straightforward using the tables \cite{24}. 

Let $(V,\omega )$ be a symplectic module of a complex semisimple group
$G$ which admits a Lagrangian highest weight vector orbit $\cal C$
and $G^0$ a connected real form. There exists a $G^0$-invariant
real structure $\tau$ on $V$ if and only if ${\fr g}^0 = Lie \, G^0$
occurs in the list of Proposition \ref{realformProp}. Consider first
the case where $V$ is irreducible. Up to multiplying $\omega$ with
$\lambda \in \mbox{\C}^{\ast}$ we can assume that $\omega$ and $\tau$ are
compatible, i.e.\ $\omega |V^{\tau}$ is a real symplectic structure.
The  nondegenerate $G^0$-invariant  Hermitian form
$\gamma = \sqrt{-1}\omega (\cdot ,\tau \cdot )$ is unique up to scaling by
$\lambda \in \mbox{\R}^{\ast}$.

Consider now the reducible case ${\fr g} = {\fr sl}_n(\mbox{\C})$ or
${\fr g} = {\fr sp}_n(\mbox{\C})$ and $V = U \oplus
U^{\ast}$, $U$ the standard $\fr g$-module. Let ${\cal C} \subset V$ 
be a Lagrangian
highest weight vector orbit. Then $L = {\cal C} \cup \{ 0\} \subset V$ is
a $G$-invariant linear subspace. Let $\tau$ be a $G^0$-invariant
real structure on $V$ compatible with $\omega$ (such pair $(\omega , \tau )$
exists). The corresponding  $G^0$-invariant Hermitian form $\gamma |L$
is unique up to scaling by $\lambda \in \mbox{\R}$. It is nondegenerate
if and only if $V = L \oplus \tau L$ and zero otherwise, i.e.\ if
$\tau L = L$. In the first (respectively second) case $V$ is the sum of two
irreducible complex or quaternionic (respectively real) $G^0$-modules.
This implies that if $\gamma |L$ is nondegenerate then ${\fr g}^0 =
{\fr su}_{k,l}$ if ${\fr g} = {\fr sl}_n(\mbox{\C})$ and ${\fr g}^0 =
{\fr sp}(k,l)$ if  ${\fr g} = {\fr sp}_n(\mbox{\C})$.

Now we return to the general case, i.e.\ $V, G, G^0, \omega , \tau ,\gamma$
and $\cal C$ are as above and the $G$-module $V$ may be reducible or
irreducible. Assume that $G^0$ has an open orbit $M = G^0 p \subset
P({\cal C})$ with compact stabilizer $H^0$. Then by Theorem \ref{redThmII}
the corresponding complex Lie group $H = (H^0)^{\mbox{\C}} \subset G$
is precisely the maximally reductive subgroup of the parabolic
stabilizer $G_p = P$. Its Lie algebra ${\fr h} = {\fr h}' \oplus {\fr t}$
is the direct sum of the semisimple Lie algebra ${\fr h}'$ (specified
in Theorem 
\ref{listthm} and Table 1) and the centre $\fr t$ which has dimension
\[ \dim {\fr t} = {\mbox{rk}} \, {\fr g} - {\mbox{rk}} \, {\fr h}' =
\left\{ \begin{array}{ll}
1 & \mbox{if}\; {\fr g}\; \mbox{is simple}\\
2 & \mbox{else.}\end{array}\right.\]
The real form ${\fr g}^0$ of $\fr g$ contains the compact real form
${\fr h}^0$ of $\fr h$. This is only the case for the real forms included
in the next proposition.

\begin{prop} \label{goodrealformsProp} Let $V$ be a $G$-module as in
Proposition  \ref{realformProp}
and $\fr h$ the maximally reductive subalgebra of ${\fr p} = Lie \, P$,
where $P = G_p$ is the stabilizer of a highest weight direction
$p\in P(V)$. Any real form ${\fr g}^0$ listed in Proposition \ref{realformProp}
which contains a subalgebra isomorphic to the compact form
${\fr h}^0$ of $\fr h$ is (up to automorphism of $\fr g$) listed in one of
the triples $(V,{\fr g}^0, {\fr h}^0)$ below. For the simple Lie algebras
$\fr g$ we have:
\begin{enumerate}
\item[A)] $(\vee^3 \mbox{\C}^2, {\fr sl}_2 (\mbox{\R}), {\fr so}_2)$,
$(\wedge^3  \mbox{\C}^6, {\fr su}_{3,3}, {\fr s}({\fr u}_3 \oplus
{\fr u}_3))$\\
and $(\mbox{\C}^n \oplus (\mbox{\C}^n)^{\ast}, {\fr su}_n$ or
${\fr su}_{1,n-1},{\fr s}({\fr u}_1 \oplus {\fr u}_{n-1}))$.
\item[B)] (Spinor module, ${\fr so}_{1,10}, {\fr u}_5)$.
\item[C)] $(V(\pi_3), {\fr sp}_3(\mbox{\R}), {\fr u}_3)$ and\\
$(\mbox{\C}^{2n}\oplus (\mbox{\C}^{2n})^{\ast} \cong 2V(\pi_1),
{\fr sp}(n)$ or ${\fr sp} (1,n-1), {\fr sp}(n-1) \oplus {\fr u}_1)$.
\item[D)] $(V(\pi_6), {\fr so}_{12}^{\ast},{\fr u}_6)$.
\item[E)] $(V(\pi_1), {\fr e}_7^{(-25)}, {\fr e}_6 \oplus {\fr so}_2)$.
\end{enumerate}
The remaining cases are:
$(\mbox{\C}^2 \otimes \mbox{\C}^n, {\fr sl}_2(\mbox{\R}) \oplus
{\fr so}_{2,n-2}, {\fr so}_2\oplus {\fr so}_2\oplus {\fr so}_{n-2})$\\
and $(\mbox{\C}^2 \otimes \mbox{\C}^7, {\fr sl}_2(\mbox{\R}) \oplus
{\fr g}_2^{(2)}, {\fr so}_2\oplus {\fr u}_2)$.
\end{prop}
  
 We will show now that each of these triples really  corresponds to a special
 pseudo-K\"ahler homogeneous manifold $M = G^0/H^0$ with compact stabilizer.
 More precisely, we prove the following main theorem.  
  
\begin{thm}\label{mainThm} Let $(V,{\fr g}^0, {\fr h}^0)$
be any of the triples listed
in Proposition \ref{goodrealformsProp},  $G$, $G^0$ and $H^0$ the connected
linear Lie groups with Lie algebras ${\fr g} = {\fr g}^0\otimes {\mbox{\C}}$,
${\fr g}^0$ and ${\fr h}^0$ respectively. Then $V$ has a Lagrangian
highest weight vector orbit ${\cal C} \subset V$ such that there exists
an open orbit of $G^0$ on $P({\cal C}) \subset P(V)$ with compact
stabilizer $H^0$ and any such orbit is a homogeneous special 
pseudo-K\"ahler manifold of the group $G^0$. 

    Conversely, any homogeneous
special pseudo-K\"ahler manifold of a real semisimple group with compact
stabilizer arises in this way. Moreover, it is homothetic to one of the
following pseudo-Hermitian locally symmetric spaces.

 The special pseudo-K\"ahler metric
is positively defined in the cases i) and ii):
\begin{enumerate}
\item[i)] ${\mbox{\C}}P^n= SU_{n+1}/S(U_1 \cdot U_n )  :\quad $ 
$({\mbox{\C}}^{n+1} \oplus (\mbox{\C}^{n+1})^{\ast}, {\fr su}_{n+1},
{\fr s}({\fr u}_1 \oplus {\fr u}_n)) .$
\item[ii)] $ {\mbox{\C}}P^{2n+1} = Sp_{n+1}/U_1 \cdot Sp_n : \quad $
 $(\mbox{\Ha}^{n+1} \oplus
(\mbox{\Ha}^{n+1})^{\ast}, {\fr sp}(n+1),{\fr sp}(n)\oplus {\fr u}_1)$.
\end{enumerate}
The cases iii)-ix) correspond to formal moduli spaces, i.e.\ the
special pseudo-K\"ahler metric is negative definite:
\begin{enumerate}
\item[iii)] $ {\mbox{\C}}H^n= SU_{1,n}/S(U_1 \cdot U_n )  : \quad $
$(\mbox{\C}^{n+1} \oplus (\mbox{\C}^{n+1})^{\ast}, {\fr su}_{1,n},
{\fr s}({\fr u}_1 \oplus {\fr u}_n))$.
\item[iv)] ${\mbox{\C}} H^1 =SL_2(\mbox{\R})/SO_2 : \quad $
 $(\vee^3 \mbox{\C}^2, {\fr sl}_2 (\mbox{\R}),
{\fr so}_2)$.
\item[v)] $SU_{3,3}/S(U_3\cdot U_3) : \quad $
$(\wedge^3  \mbox{\C}^6, {\fr su}_{3,3}, {\fr s}({\fr u}_3 \oplus
{\fr u}_3))$.
\item[vi)] $Sp_3(\mbox{\R})/U_3 : \quad $
$(V(\pi_3), {\fr sp}_3(\mbox{\R}), {\fr u}_3)$.
\item[vii)] $SO_{12}^{\ast}/U_6 : \quad $
$(V(\pi_6), {\fr so}_{12}^{\ast},{\fr u}_6)$.
\item[viii)] $E_7^{(-25)}/E_6\cdot SO_2 : \quad $
$(V(\pi_1), {\fr e}_7^{(-25)}, {\fr e}_6 \oplus {\fr so}_2)$.
\item[ix)] ${\mbox{\C}}H^1 \times (SO_{2,n} /SO_2 \cdot SO_n) : 
\quad$ 
$(\mbox{\C}^2 \otimes \mbox{\C}^{n+2}, {\fr sl}_2(\mbox{\R}) \oplus
{\fr so}_{2,n}, {\fr so}_2\oplus {\fr so}_2\oplus {\fr so}_{n})$.
\end{enumerate}
In the remaining cases the metric is indefinite of (complex) signature (10,5),
(1,2n) and (1,5) in the cases x), xi) and xii) respectively.
\begin{enumerate}
\item[x)] $Spin_{1,10}/U_5 \subset SO_{2,10}/U_{1,5} : \quad$ (spinor module,
${\fr so}_{1,10}, {\fr u}_5)$.
\item[xi)] $Sp(1,n)/U_1\cdot Sp(n) \subset 
SU_{2,2n}/S(U_1 \cdot U_{1,2n}) : \quad$
$(\mbox{\Ha}^{n+1} \oplus
(\mbox{\Ha}^{n+1})^{\ast}, {\fr sp}(1,n), {\fr sp}(n) \oplus {\fr u}_1)$.
\item[xii)] $(SL_2(\mbox{\R})/SO_2) \times (G_2^{(2)}/SO_2 \cdot U_2)
\subset (SL_2(\mbox{\R})/SO_2) \times (SO_{3,4} /SO_2 \cdot SO_{1,4}):
\quad  (\mbox{\C}^2 \otimes \mbox{\C}^7, {\fr sl}_2(\mbox{\R}) \oplus
{\fr g}_2^{(2)}, {\fr so}_2\oplus {\fr u}_2)$.
\end{enumerate}
\end{thm}

\noindent
{\bf Proof:} First we treat the case of reducible modules: ${\fr g} =
{\fr g}^0 \otimes \mbox{\C} = {\fr sl}_{n+1} (\mbox{\C})$ or
${\fr g} = {\fr sp}_{n+1}(\mbox{\C})$ and the symplectic $\fr g$-module
$(V,\omega )$ admits the $\fr g$-invariant Lagrangian splitting
$V = L \oplus L'$, where $L$ is isomorphic to the standard $\fr g$-module
and $L'$ is isomorphic to its dual. Let $\tau$ be a real structure on $V$
invariant under the real form ${\fr g}^0 = {\fr su}_{k,l}$
(respectively ${\fr sp}_{k,l}$) of ${\fr g} = {\fr sl}_{n+1} (\mbox{\C})$
(respectively ${\fr g} = {\fr sp}_{n+1}(\mbox{\C})$), $k+l = n+1$.
As discussed above, we can choose the data $(\omega , L, L',\tau )$ such that
$\tau$ is compatible with $\omega$ and interchanges $L$ and $L'$, i.e.\
$L' = \tau L$.  The corresponding ${\fr g}^0$-invariant nondegenerate
Hermitian structure $\gamma |L$,
$\gamma = \sqrt{-1}\omega (\cdot ,\tau \cdot )$, is unique up to scaling by
$\lambda \in \mbox{\R}^{\ast}$. Up to multiplying $\tau$ by $-1$
we can assume that $\gamma |L$   has signature $(k,l)$ (respectively
$(2k,2l)$) if ${\fr g}^0 = {\fr su}_{k,l}$ (respectively
${\fr sp}(k,l)$).

Let $G^0$ be the linear group with Lie algebra
${\fr g}^0$. We consider its action on the projectivization $P({\cal C})$
of the highest weight vector orbit ${\cal C} = L- \{ 0\}$.
The open $G^0$-orbits $M = G^0p$ on $P({\cal C})$ are precisely the orbits
of non-null directions $p = \mbox{\C}^{\ast}v$, i.e.\ $\gamma (v,v) \neq 0$.
The cone ${\cal C}_M \subset V$ over any open orbit $M$ is a special cone,
see Def.~\ref{LagrangianDef}. The pseudo-K\"ahler metric $\gamma |L$
is invariant under the group $SU_{k,l}$ if the signature of
$\gamma |L$ is $(k,l)$. This shows that, up to a scaling of the metric by
a {\em positive} factor, the special pseudo-K\"ahler manifold
$(M = G^0p,g)$ is isometric to the pseudo-Hermitian symmetric space
$SU_{k,l}/S(U_1 \cdot U_{k-1,l})$ (respectively $SU_{k,l}/S(U_{k,l-1}
\cdot U_1)$) of signature $(k-1,l)$ (respectively $(l-1,k)$) if
$\gamma (v,v) >0$ (respectively $\gamma (v,v) <0$), $p =\mbox{\C}^{\ast}v$.

It follows that if $G^0 = SU_{k,l}$ and $H^0$ is the stabilizer of $p$
in $G^0$ then $H^0$ is compact if and only if the special 
pseudo-K\"ahler metric $g$ on $M = G^0p \subset P({\cal C})$ is of definite
signature. Moreover, $(M,g)$ is homothetic to complex projective space
${\mbox{\C}}P^n$ (respectively complex hyperbolic space ${\mbox{\C}}H^n$)
if $g>0$ (respectively $g<0$). Up to an automorphism of ${\fr g} =
{\fr sl}_{n+1} (\mbox{\C})$ we obtain the cases i) and iii) of the
theorem.

If $G^0 = Sp(k,l)$ then $H^0$ is compact if and only if either $G^0$
is compact, the special pseudo-K\"ahler metric $g$ on 
$M = G^0p \cong Sp(n+1)/Sp(n)\cdot U_1$ is positive
definite and $(M,g)$ is homothetic to ${\mbox{\C}}P^{2n+1} = SU_{2n+2}/
S(U_{2n+1} \cdot U_1)$ or $G^0 \cong Sp(1,n)$, the special 
pseudo-K\"ahler metric on $M = G^0p \cong Sp(1,n)/Sp(n)\cdot U_1$ is negative
definite and $(M,g)$ is homothetic to an open $Sp(1,n)$-orbit on 
the pseudo-Hermitian symmetric
space $SU_{2,2n}/S(U_1 \cdot U_{1,2n})$.  This covers the cases ii)
and xi) and completes the discussion of the reducible modules.

Next we study the irreducible modules. So let $(V,{\fr g}^0,{\fr h}^0)$ be
a triple listed in Proposition 
\ref{goodrealformsProp} such that ${\fr g}^0$
is a noncompact real form of the complex semisimple Lie algebra $\fr g$
and $V$ is an irreducible symplectic $\fr g$-module. The highest weight vector
orbit ${\cal C} = Gv \subset V$ of the complex linear group $G$ with
$Lie \, G = {\fr g}$ is a Lagrangian cone. The stabilizer $G_p = P$ of
$p = \mbox{\C}^{\ast}v \in P(V)$ in $G$ is a parabolic subgroup with maximally
reductive subgroup $H$. The maximal rank subalgebra ${\fr h}^0 \subset
{\fr g}^0$ is the compact real form of ${\fr h} = Lie \, H$.
Any Cartan subalgebra ${\fr a}^0 \subset {\fr h}^0$ is a Cartan subalgebra
of ${\fr g}^0$   and its complexification ${\fr a} = {\fr a}^0
\otimes \mbox{\C} \subset {\fr h} \subset {\fr g}$ is a  Cartan subalgebra
of $\fr g$. Let $R^+$ be a system of positive roots of $\fr g$
with respect to $\fr a$ and $p$ the highest weight direction of the
irreducible $\fr g$-module $V$ with respect to $R^+$. We can write
\[ {\fr h} = {\fr a} + \sum_{\alpha \in R_0} {\fr g}_{\alpha} \]
for some root subsystem $R_0 \subset R = R^+ \cup R^-$ and
\[ {\fr p} = Lie \, P =  {\fr h} +
\sum_{\alpha \in R^+ - R_0}{\fr g}_{\alpha}\, ,
\quad R^+_0 = R_0 \cap R^+\, . \]
We denote by $\sigma$ the complex antilinear involution of $\fr g$ with
fixed point set ${\fr g}^{\sigma} = {\fr g}^0$. Since ${\fr a}^0
\subset {\fr g}^0$ is a compact Cartan subalgebra, we have that
\[ {\fr a}^{\mbox{\R}} := \{ \alpha \in R | \alpha ({\fr a}) \subset
\mbox{\R} \} = i{\fr a}^0\, ,\quad \sigma {\fr g}_{\alpha} =
{\fr g}_{-\alpha}\]
and hence ${\fr p} \cap \sigma {\fr p} = {\fr h}$ and  ${\fr p} \cap {\fr g}^0
= {\fr h}^0$. This shows that $M = G^0p \cong G^0/H^0$ is open in
$\overline{M} = Gp = P({\cal C})$.  Moreover, any $G^0$-orbit on
$P({\cal C})$ with stabilizer $H^0$ is obtained in this way. In fact, let
$\fr p$ be the Lie algebra of the stabilizer in $G$ of a point
$p\in P({\cal C})$ and assume that ${\fr p}\cap {\fr g}^0 = {\fr h}^0$.
The point $p$ is a highest weight direction for some choice of Cartan
subalgebra ${\fr a} \subset {\fr h}$ and system of positive roots
$R^+$ of $\fr g$. Since  any two Cartan
subalgebras of ${\fr h}$ are conjugated by an element of $H$, we can
assume that ${\fr a} = {\fr a}^0 + i{\fr a}^0$, ${\fr a}^0 = {\fr a} \cap
{\fr h}^0$.

The next lemma shows that, under our assumptions, the cone
${\cal C}_M \subset V$ over $M$ is special and the signature of
the Hermitian metric $\gamma$ on ${\cal C}_M$ is easily computed.

\begin{lemma} \label{signLemma} Let $(V,\omega)$ be an irreducible symplectic
$G$-module of a complex semisimple Lie group $G$ with Lagrangian highest
weight vector orbit ${\cal C} \subset V$ and $G^0$ a connected real
form of $G$ for which there exists a $G^0$-invariant real structure
$\tau$ on $V$. Assume that the stabilizer $H^0$ of $p\in P({\cal C})$ 
in $G^0$ is a compact real form of the maximally reductive
subgroup $H$ of the stabilizer $G_p = P$ of $p$ in $G$.
Then the orbit $M = G^0 p \subset P({\cal C})$ is open
and the cone ${\cal C}_M \subset V$ over $M$ is special in the sense
of Def.~\ref{LagrangianDef}. The Hermitian form
$\gamma = \sqrt{-1}\omega (\cdot ,\tau \cdot )$  has signature
(k,l) or (l,k) on ${\cal C}_M$, where k (respectively l) is the number
of compact (respectively noncompact) roots in $R^+ - R^+_0$. Here
$R^+$ (respectively $R_0^+ \subset R^+$) is a system of positive
roots of $\fr g$ (respectively $\fr h$) with respect to the Cartan
subalgebra ${\fr a} = {\fr a}^0 \otimes \mbox{\C}$, ${\fr a}^0 \subset
{\fr h}^0$ a Cartan subalgebra of ${\fr g}^0$. If $\gamma$ has signature
$(k,l)$ on ${\cal C}_M$, then the special pseudo-K\"ahler metric
$g$ of $M = G^0 p$, $p = \mbox{\C}^{\ast} v$, has signature $(k-1,l)$
if $\gamma (v,v) > 0$ and signature $(l-1,k)$ if $\gamma (v,v) < 0$.
\end{lemma}

\noindent
{\bf Proof:} As above, we denote by $\sigma$ the {\C}-antilinear involution
of $\fr g$ with fixed point set  ${\fr g}^{\sigma} = {\fr g}^0$ and
recall that $\sigma {\fr g}_{\alpha} =
{\fr g}_{-\alpha}$ for all $\alpha \in R$, where $R$ is the system
of roots of $\fr g$ with respect to $\fr a$. We choose a system of positive
roots $R^+\subset R$ such that our base point $v\in M$ is a highest weight
vector with respect to $R^+$ and denote by $\Lambda$ the  highest weight.
Then ${\fr h} = {\fr h}^0 \otimes \mbox{\C} = {\fr a} + \sum_{\alpha \in
R_0}{\fr g}_{\alpha}$ for some root subsystem $R_0 = R_0^+ \cup -R_0^+
\subset R$, $R_0^+ = R^+ \cap R_0$.
Then the tangent space
$T_v{\cal C}_M$ is precisely
\[ T_v {\cal C}_M = {\fr g}v = \mbox{\C}v \oplus
\bigoplus_{\alpha \in R^+ - R_0^+}
{\fr g}_{-\alpha}v = V_{\Lambda} \oplus \bigoplus_{\alpha \in R^+ - R_0^+}
V_{\Lambda -\alpha}\, .\]
The ${\fr g}^0$-invariance of the real structure $\tau$ on $V$
is equivalent to $\tau (xu) = \sigma (x) \tau (u)$ for all
$x\in {\fr g}$, $u\in V$ and hence we have that:
\[ \tau T_v {\cal C}_M = \mbox{\C}\tau v  \oplus \bigoplus_{\alpha \in R^+ - R_0^+}
{\fr g}_{\alpha}\tau v = V_{-\Lambda} \oplus
\bigoplus_{\alpha \in R^+ - R_0^+}
V_{\alpha -\Lambda}\, .\]
This shows that $V = T_v {\cal C}_M \oplus \tau T_v {\cal C}_M$ and therefore
$\gamma$ is nondegenerate on $L = T_v {\cal C}_M$.
Moreover, any $G^0$-invariant Hermitian form on $L$ is
proportional to $\gamma |L$ over the reals.  A straightforward calculation
shows that $\gamma |L$
has the orthogonal basis
$(v, e_{\alpha} v, \alpha \in R^- - R_0^-)$, $R^- = - R^+$, $R_0^- = -R_0^-$,
which satisfies:
\[ \gamma (e_{\alpha} v,e_{\alpha} v) =
\Lambda ([ e_{\alpha},\sigma e_{\alpha}]) \gamma (v,v) =
\langle \alpha , \Lambda \rangle \langle e_{\alpha},\sigma e_{\alpha} \rangle
\; \left\{  \begin{array}{ll}
>0& \mbox{if}\; \alpha \; \mbox{is compact}\\
<0& \mbox{if}\; \alpha \; \mbox{is noncompact.}\end{array}\right.\]
Here we have used that the Killing form $\langle \cdot , \cdot \rangle$
is negative (respectively positive) definite on
$({\fr g}^{\alpha} + {\fr g}^{-\alpha})^{\sigma}$ if $\alpha$ is a
compact (respectively noncompact) root and the fact that $R^- - R_0^- =
\{ \alpha \in R| \langle \alpha , \Lambda \rangle < 0 \}$. This
proves the lemma. 

Now we complete the proof of the theorem. Thanks to Lemma \ref{signLemma}
it only remains to show that the special pseudo-K\"ahler metric on $M
= G^0 p \cong G^0/H^0$ is pseudo-Hermitian locally symmetric. This is immediate
if the isotropy representation is irreducible. In the remaining
two cases we have the following inclusions $({\fr g}^0, {\fr h}^0)
\subset (\tilde{{\fr g}}^0,\tilde{{\fr h}}^0)$ of pairs:
\[ ({\fr so}_{1,10}, {\fr u}_5) \subset ({\fr so}_{2,10}, {\fr u}_{1,5})
\quad \mbox{and} \]
\[ ({\fr sl}_2(\mbox{\R}) \oplus {\fr g}_2^{(2)},{\fr so}_2 \oplus {\fr u}_2)
\subset ({\fr sl}_2(\mbox{\R}) \oplus {\fr so}_{3,4}, {\fr so}_2 \oplus
{\fr so}_2 \oplus {\fr so}_{1,4})\, .\]
The irreducible symplectic $\fr g$-module $V$ ($V =$ spinor module
in the case ${\fr g} = {\fr so}_{11}(\mbox{\C})$ and
$V = \mbox{\C}^2 \otimes  \mbox{\C}^7$ for
${\fr g} = {\fr sl}_2(\mbox{\C}) \oplus {\fr g}_2(\mbox{\C})$)
carries the structure
of irreducible symplectic $\tilde{{\fr g}}$-module, $\tilde{{\fr g}} =
\tilde{{\fr g}}^0 \otimes \mbox{\C}$ (semispinor module if $\tilde{{\fr g}}
= {\fr so}_{12}(\mbox{\C})$, tensor product of standard modules if
$\tilde{{\fr g}}
= {\fr sl}_2(\mbox{\C}) \oplus {\fr so}_7(\mbox{\C})$), and a
$\tilde{{\fr g}}^0$-invariant
real structure, as follows from Proposition \ref{realformProp}.
This shows that the special pseudo-K\"ahler metric of $M = G^0 p \cong
G^0/H^0 \hookrightarrow \tilde{G}^0 /\tilde{H}^0$ is canonically
extended to a $\tilde{G}^0$-invariant special pseudo-K\"ahler metric
on $\tilde{M} = \tilde{G}^{0}/\tilde{H}^{0}$. Here $\tilde{G}^0$ 
and $\tilde{H}^0$ are the connected
linear groups with $Lie \, G^0 = {\fr g}^0$ and $Lie \, H^0
= {\fr h}^0$. On the other hand, $\tilde{M} = 
\tilde{G}^0 /\tilde{H}^0$ carries the
structure of an 
{\em irreducible} pseudo-Hermitian symmetric space. Hence 
the special pseudo-K\"ahler metric of $\tilde{M}$ is 
pseudo-Hermitian symmetric. Since the special pseudo-K\"ahler manifold 
$M$ is isometrically embedded into 
$\tilde{M}$ as an open $G^0$-orbit it follows that $M$ is 
pseudo-Hermitian locally symmetric 
finishing the proof of Theorem~\ref{mainThm}.


\begin{thebibliography}{10}
\bibitem[1]{1}  D.V.\ Alekseevsky:  Classification of
quaternionic spaces with a transitive solvable group of motions, {\it 
Math.\ USSR Izvestija} {\bf 9}, No.\ 2 (1975), 297-339.
\bibitem[2]{2}  D.V.\ Alekseevsky: Flag manifolds,
{\it Zbornik Radova} {\bf  6} (14) (1997), 3-35;
Preprint ESI, No.\ 415 (1997), 32p.
\bibitem[3]{3} D.V.\ Alekseevsky, V.\ Cort\'es:   Isometry groups of
homogeneous quaternionic K\"ahler manifolds, to appear in
{\it Journal of Geometric Analysis}; available as preprint Erwin
Schr\"odinger Institut  230 (1995).  
\bibitem[4]{4} W.M.\
Boothby: Homogeneous complex contact manifolds,
{\it Proceedings of the Symposia in Pure Mathematics} {\bf 3} (1961),
144-154.
\bibitem[5]{5} W.M.\ Boothby:  A note on
homogeneous complex contact manifolds, {\it Proc.\ Amer.\ Math.\
Soc.} {\bf 13} (1962), 276-280.
\bibitem[6]{6} S.\ Cecotti:
Homogeneous K\"ahler manifolds and T-algebras in 
$N = 2$ supergravity and superstrings, {\it Commun.\ Math.\ Phys.} {\bf
124} (1989), 23-55.  
\bibitem[7]{7} V.\ Cort\'es: Alekseevskian spaces, {\it Diff.\ 
Geom.\ Appl.} {\bf 6} (1996),
129-168.
\bibitem[8]{8} V.\ Cort\'es: Homogeneous
special geometry, {\it Transformation Groups} {\bf 1}, No.\ 4 (1996),
337-373.  
\bibitem[9]{9} V.\ Cort\'es:  On hyper K\"ahler
manifolds associated to Lagrangian K\"ah\-ler submanifolds of 
$T^*\mbox{\C}^n$, {\it Trans.\ Amer.\ Math.\ Soc.} {\bf 350} (1998), 
no.\ 8, 3193-3205. 
\bibitem[10]{10} B.\ de Wit, A.\ Van Proeyen:  Potentials
and symmetries of general gauged $N = 2$ supergravity-Yang-Mills
models, {\it Nucl.\ Phys.} {\bf B245} (1989), 89-117.
\bibitem[11]{11}  B.\ de Wit, A.\ Van Proeyen: Special
geometry, cubic polynomials and homogeneous quaternionic spaces,
{\it Commun.\ Math.\ Phys.} {\bf 149} (1992), 307-333.
\bibitem[12]{12} B.\ de Wit, F.\ Vanderseypen, A.\ Van Proeyen:
Symmetry structure of special
geometries, {\it Nucl.\ Phys.} {\bf B400} (1993), 463-521.
\bibitem[13]{13} D.S.\ Freed:  Special K\"ahler manifolds, 
hep-th/9712042
\bibitem[14]{14} S.\ Helgason:  Differential geometry, Lie groups and
symmetric spaces, Academic Press, Orlando, 1978.
\bibitem[15]{15} N.J.\ Hitchin: The moduli space of complex
Lagrangian submanifolds, preprint math.DG/9901069. 
\bibitem[16]{16} J.-M.\ Hwang:  Rigidity of homogeneous contact
manifolds under Fano deformation (preprint).
\bibitem[17]{17} E.\ Looijenga, V.A.\  Lunts : A Lie algebra 
attached to a projective variety, {\it Invent.\ Math.} 
{\bf 129} (1997), 361-412.
\bibitem[18]{18} Zhiqin Lu:  A note on special K\"ahler manifolds,
Columbia University preprint (1998). 
\bibitem[19]{19} S.A.\ Merkulov, L.J.\ Schwachh\"ofer: 
Twistor solution of the holonomy problem, in Geometric Issues
in the Foundation of Science, Oxford Univ.\ Press, 1998.  
\bibitem[20]{20} A.L.\ Onishchik: Topology of transitive transformation
groups, Barth Verlagsgesellschaft, Leipzig, Berlin, Heidelberg, 1994.
\bibitem[21]{21} A.L.\ Onishchik, E.B.\ Vinberg:  Lie groups and
algebraic groups,
Springer, Berlin, Heidelberg, 1990.
\bibitem[22]{22} A.L.\ Onishchik, E.B.\ Vinberg (Eds.):
Lie groups and Lie algebras III (EMS v.\ 41), Springer,
Berlin, Heidelberg, New York, 1994.
\bibitem[23]{23} L.J.\ Schwachh\"ofer: 
On the classification of holonomy representations, Habilitationsschrift,
Leipzig, 1998.
\bibitem[24]{24} J.\ Tits:  Tabellen zu den einfachen Lieschen Gruppen
und ihren Darstellungen, LNM 40, Springer, Berlin, 1967.
\bibitem[25]{25} J.A.\ Wolf:  Complex homogeneous contact manifolds
and quaternionic symmetric spaces, {\it J.\ Math.\ Mech.} {\bf 14}
(1965), 1033-1047.
\bibitem[26]{26} J.A.\ Wolf: The action of a real semisimple group
on a complex flag manifold. I: Orbit structure and holomorphic arc
components, {\it Bulletin of the AMS} {\bf 75}, No.\ 6 (1969), 1121-1237.
\bibitem[27]{27} F.L.\ Zak: Tangents and secants of algebraic varieties
(Translations of mathematical monographs, v.\ 127), AMS, Providence, 1993.
\end{thebibliography}
\end{document}